\documentclass[12pt,a4paper]{article}
\usepackage[width=18cm,height=23cm]{geometry}
\usepackage{amsmath,amssymb,amsfonts,amsthm,enumerate,float}

\usepackage{bbm}
\usepackage{graphicx}
\usepackage{tikz}
\usetikzlibrary{shapes,backgrounds}
\usepackage{tikz-cd}
\usetikzlibrary{%
  matrix,%
  calc,%
  arrows%
}

\usepackage{mathrsfs,amsbsy,bm,mathtools}
\usepackage{indent first}

\newcommand{\bc}{\mathbb{C}}
\newcommand{\bd}{\mathbb{D}}
\newcommand{\cc}{\mathscr{C}}
\newcommand{\cl}{\mathscr{L}}
\newcommand{\bp}{\mathbb{P}}
\newcommand{\bu}{\mathbb{U}}
\newcommand{\bh}{\mathbb{H}}
\newcommand{\bb}{\mathbb{B}}
\newcommand{\bt}{\mathbb{T}}
\newcommand{\bz}{\mathbb{Z}}

\newcommand{\br}{\mathbb{R}}
\newcommand{\bq}{\mathbb{Q}}

\newcommand{\cf}{\mathscr{F}}

\newcommand{\cond}[1]{\quad{\scriptstyle{(#1)}}}
\newcommand{\rw}{\rightarrow}
\newcommand{\lrw}{\longrightarrow}

\newtheorem{Thm}{Theorem}[section]%
\newtheorem{Lem}[Thm]{Lemma}%
\newtheorem{Cor}[Thm]{Corollary}%
\newtheorem{Prop}[Thm]{Proposition}%

\theoremstyle{definition} 
\newtheorem{Def}[Thm]{Definition}%

\title{\bf Directed harmonic currents near\\
 non-hyperbolic linearized singularities}
\author{Zhangchi Chen}
\date{\today}

\begin{document}
\maketitle

\abstract{
Let $(\bd^2,\cf,\{0\})$ be a singular holomorphic foliation on the unit bidisc $\bd^2$ defined by the linear vector field
\[
z \,\frac{\partial}{\partial z}+ \lambda \,w \,\frac{\partial}{\partial w},
\]
where $\lambda\in\bc^*$. Such a foliation has a non-degenerate linearized singularity at $0$. Let $T$ be a harmonic current directed by $\cf$ which does not give mass to any of the two separatrices $(z=0)$ and $(w=0)$ and whose the trivial extension $\tilde{T}$ across $0$ is $dd^c$-closed. The Lelong number of $T$ at $0$ describes the mass distribution on the foliated space. In 2014 Nguyen proved that when $\lambda\notin\br$, i.e. $0$ is a hyperbolic singularity, the Lelong number at $0$ vanishes. For the non-hyperbolic case $\lambda\in\br^*$ the article proves the following results. The Lelong number at $0$:

1) is strictly positive if $\lambda>0$;

2) vanishes if $\lambda\in\mathbb{Q}_{<0}$;

3) vanishes if $\lambda<0$ and $T$ is invariant under the action of some cofinite subgroup of the monodromy group.}


%%%%%%%%%%%%%%%%%%%%%%%%%%%%%
\section{Introduction}
\label{sect-intro}
%%%%%%%%%%%%%%%%%%%%%%%%%%%%%
The dynamical properties of holomorphic foliations have drawn great attention recently \cite{Nguyen-2020-05}. Let us see one of the recent interesting results.

\begin{Thm}[Dinh-Nguy\^en-Sibony \cite{Dinh-Nguyen-Sibony-2018}] Let $\cf$ be a holomorphic foliation with only hyperbolic singularities in a compact K\"ahler surface $(X,\omega)$. Assume that $\cf$ admits no directed positive closed current. Then there exists a unique positive $dd^c$-closed current $T$ of mass $1$ directed by $\cf$.
\end{Thm}

The first version was stated for $X=\bp^2$ and proved by Forn\ae ss-Sibony \cite{Fornaess-Sibony-2010}. Later Dinh-Sibony proved unique ergodicity for foliations in $\bp^2$ with an invariant curve \cite{Dinh-Sibony-2018}. So one may expect to describe recurrence properties of leaves by studying the density distribution of directed harmonic currents. One has the following the result about leaves.

\begin{Thm}[Forn\ae ss-Sibony \cite{Fornaess-Sibony-2010}]\label{FS-Mass} Let $(X,\cf,E)$ be a hyperbolic foliation on a compact complex surface $X$ with singular set $E$. Assume that
\begin{enumerate}
\item there is no invariant analytic curve;
\item all the singularities are hyperbolic;
\item there is no non-constant holomorphic map $\bc\rw X$ such that out of $E$ the image of $\bc$ is locally contained in a leaf.
\end{enumerate}
Then every harmonic current $T$ directed by $\cf$ gives no mass to each single leaf.
\end{Thm}

A practical way to measure the density of harmonic currents is to use the notion of Lelong number introduced by Skoda \cite{Skoda-1982}. Indeed Theorem~\ref{FS-Mass} above is equivalent to the statement that the Lelong number of $T$ vanishes everywhere outside $E$. Another result holds near hyperbolic singularities. 

\begin{Thm}[Nguy\^en \cite{Nguyen-2018}] Let $(\bd^2,\cf,\{0\})$ be a holomorphic foliation on the unit bidisc $\bd^2$ defined by the linear vector field $Z(z,w)=z\,\frac{\partial}{\partial z}+\lambda\,w\,\frac{\partial}{\partial w},$ where $\lambda\in\bc\backslash\br$, that is to say, $0$ is a hyperbolic singularity. Let $T$ be a harmonic current directed by $\cf$ which does not give mass to any of the two separatrices $(z=0)$ and $(w=0)$. Then the Lelong number of $T$ at $0$ vanishes.
\end{Thm}

Nguy\^en proved that the Lelong number of any directed harmonic current which gives no mass to invariant hyperplanes, vanishes near {\sl weakly hyperbolic} singularities in $\bc^n$ \cite{Nguyen-2020-05}. This result is optimal, see \cite{Dinh-Wu-2020}. The mass-distribution problem would be completed once the behaviour of harmonic currents near non-hyperbolic and near degenerate singularities would be understood.

The present paper answers (partly) the problem in the non-hyperbolic linearizable singularity case.

\begin{Thm} Let $(\bd^2,\cf,\{0\})$ be a holomorphic foliation on the unit bidisc $\bd^2$ defined by the linear vector field $Z(z,w)=z\,\frac{\partial}{\partial z}+\lambda\,w\,\frac{\partial}{\partial w}$, where $\lambda\in\br^*$. Let $\cf$ be a harmonic current directed by $\cf$ which does not give mass to any of the two separatrices $(z=0)$ and $(w=0)$. Then the Lelong number of $T$ at $0$
\begin{itemize}
\item is strictly positive if $\lambda>0$,
\item vanishes if $\lambda\in\bq_{<0}$.
\end{itemize}
\end{Thm}

For the concerned foliation $(\bd^2,\cf,\{0\})$, a local leaf $P_\alpha$, with $\alpha\in\bc^*$, can be parametrized by $(z,w)=(e^{-v+iu},\alpha\,e^{-\lambda v+i\lambda u})$. The {\sl monodromy group} around the singularity is generated by $(z,w)\mapsto(z,e^{2\pi i\lambda}w)$. It is a cyclic group of finite order when $\lambda\in\bq^*$, of infinite order when $\lambda\notin\bq$.

It is now ready to introduce the notion of {\sl periodic current}, an essential tool of this paper. A directed harmonic current $T$ is called {\sl periodic} if it is invariant under some cofinite subgroup of the monodromy group, i.e. under the action of $(z,w)\mapsto(z,e^{2k\pi i \lambda}w)$ for some $k\in\bz_{>0}$. Observe that if $\lambda\in\bq^*$ then any directed harmonic current is periodic. But when $\lambda\notin\bq^*$, the periodicity is a nontrivial assumption.

\begin{Thm}\label{thm:negative-periodic} Using the same notation as above, the Lelong number of $T$ at the singularity is 0 when $\lambda<0$ and the current is periodic, in particular, when $\lambda\in\bq_{<0}$.
\end{Thm}

It remains open to determine the possible Lelong number values of non-periodic $T$ when $\lambda<0$ is irrational.

%%%%%%%%%%%%%%%%%%%%%%%%%%%%%
\section{Background}
\label{sect-background}
%%%%%%%%%%%%%%%%%%%%%%%%%%%%%

To start with, recall the definition of singular holomorphic foliation on a complex surface $M$.

\begin{Def}
Let $E\subset M$ be some closed subset, possibly empty, such that $\overline{M\backslash E}=M$. A {\sl singular holomorphic foliation} $(M,E,\cf)$ consists of a holomorphic {\sl atlas} $\{(\bu_i,\Phi_i)\}_{i\in I}$ on $M\backslash E$ which satisfies the following conditions.
\begin{enumerate}[(1)]
\item For each $i\in I$, $\Phi_i: \bu_i\rw \bb_i\times\bt_i$ is a biholomorphism, where $\bb_i$ and $\bt_i$ are domains in $\bc$.
\item For each pair $(\bu_i,\Phi_i)$ and $(\bu_j,\Phi_j)$ with $\bu_i\cap\bu_j\neq\emptyset$, the transition map
\[
\Phi_{ij}:=\Phi_i\circ\Phi_j^{-1}:\Phi_j(\bu_i\cap\bu_j)\rw\Phi_i(\bu_i\cap\bu_j)
\]
has the form
\[
\Phi_{ij}(b,t)=\big(\Omega(b,t),\Lambda(t)\big),
\]
where $(b,t)$ are the coordinates on $\bb_j\times\bt_j$, and the functions $\Omega$, $\Lambda$ are holomorphic, with $\Lambda$ independent of $b$.
\end{enumerate}
\end{Def}

Each open set $\bu_i$ is called a {\sl flow box}. For each $c\in\bt_i$, the Riemann surface $\Phi_i^{-1}\{t=c\}$ in $\bu_i$ is called a {\sl plaque}. The property (2) above insures that in the intersection of two flow boxes, plaques are mapped to plaques.

A {\sl leaf} $L$ is a minimal connected subset of $M$ such that if $L$ intersects a plaque, it contains that plaque. A {\sl transversal} is a Riemann surface immersed in $M$ which is transverse to each leaf of $M$.

The local theory of singular holomorphic foliations is closely related to holomorphic vector fields. One recalls some basic concepts in $\bc^2$ (\cite{Nguyen-2020-05}, \cite{Brunella-2015}).

\begin{Def} Let $Z=P(z,w)\frac{\partial}{\partial z}+Q(z,w)\frac{\partial}{\partial w}$ be a holomorphic vector field defined in a neighborhood $\bu$ of $(0,0)\in\bc^2$. One says that $Z$ is
\begin{enumerate}[(1)]
\item {\sl singular} at $(0,0)$ if $P(0,0)=Q(0,0)=0$,
\item {\sl linear} if it can be written as
\[
Z=\lambda_1 z\frac{\partial}{\partial z}+\lambda_2 w \frac{\partial}{\partial w}
\]
where $\lambda_1$, $\lambda_2\in\bc$ are not simultaneously zero.
\item {\sl linearizable} if it is linear after a biholomorphic change of coordinates.
\end{enumerate}
\end{Def}

Suppose the holomorphic vector field $Z=P\frac{\partial}{\partial z}+Q\frac{\partial}{\partial w}$ admits a singularity at the origin. Let $\lambda_1$, $\lambda_2$ be the eigenvalues of the Jacobian matrix $\textstyle\left(\!
\begin{array}{cc}
P_z & P_w
\\
Q_z & Q_w
\end{array}
\!\right)$ at the origin.

\begin{Def} The singularity is {\sl non-degenerate} if both $\lambda_1$, $\lambda_2$ are non-zero. This condition is biholomorphically invariant.
\end{Def}

In this article, all singularities are assumed to be non-degenerate. Then the foliation defined by integral curves of $Z$ has an isolated singularity at $0$. Degenerate singularities are studied in \cite{Brunella-2015}. Seidenberg's reduction theorem \cite{Seidenberg-1968} shows that degenerate singularites can be resolved into non-degenerate ones after finitely many blow-ups.

\begin{Def}
A singularity of $Z$ is {\sl hyperbolic} if the quotient $\lambda:=\frac{\lambda_1}{\lambda_2}\in\bc^*\backslash\br$. It is {\sl non-hyperbolic} if $\lambda\in\br^*$. It is in the {\sl Poincar\'e domain} if $\lambda\notin\br_{\leqslant0}$. It is in the {\sl Siegel domain} if $\lambda\in\br_{<0}$.
\end{Def}

One can verify that the quotient is unchanged by multiplication of $Z$ by any non-vanishing holomorphic function.

One could consider $\lambda^{-1}=\frac{\lambda_2}{\lambda_1}$ instead of $\lambda$, but then $\lambda\notin\br$ iff $\lambda^{-1}\notin\br$. Thus the notion of hyperbolicity is well-defined. Also, being non-hyperbolic, in Poincar\'e domain or in Siegel domain, is well-defined. The complex number $\lambda$ will be called {\sl eigenvalue} of $Z$ at the singularity, with an inessential abuse due to this exchange $\lambda\leftrightarrow\lambda^{-1}$. The unordered pair $\{\lambda,\lambda^{-1}\}$ is invariant under local biholomorphic changes of coordinates.

Consider a holomorphic foliation $(M,E,\cf)$ where $E$ is discrete. When one tries to linearize a vector field near an isolated singularity, one has to divide power series coefficients by quantities $\lambda_1\,m_1+\lambda_2\,m_2-\lambda_j$ for $j=1,2$ where $m_1$, $m_2\in\bz_{\geqslant 1}$. To ensure convergence, these quantities have to be nonzero and not too close to zero.

The {\sl resonances} of $(\lambda_1,\lambda_2)\in\bc^2$ are defined by
\[
\mathcal{R}:=\{(m_1,m_2,j)\in\bz^3~|~m_1,\,m_2\in\bz_{\geqslant 1},~ \lambda_1\,m_1+\lambda_2\,m_2-\lambda_j=0\}.
\]
Notice that the set $\{\lambda_1\,m_1+\lambda_2\,m_2-\lambda_j~|~m_1,m_2\in\bz_{\geqslant1 }\}$ has zero as a limit point if and only if the singularity is in the Siegel domain.

We are now ready to state some linearization results in $\bc^2$.
\begin{Thm}[Poincar\'e \cite{Arnold-1988}] A singular holomorphic vector field with a non-resonant linear part, i.e. $\mathcal{R}$ is empty, such that the eigenvalue $\lambda$ is in the Poincar\'e domain, is holomorphically equivalent to its linear part.
\end{Thm}

That is to say, a singularity is linearizable if $\lambda\notin\br_{\leqslant0}\cup\bq$. To get linearization for $\lambda$ in the Siegel domain, the following result assumes the more advanced {\sl Brjuno condition}.

\begin{Thm}[Brjuno \cite{Arnold-1988,Brjuno-1979}]
A singular holomorphic vector field with a non-resonant linear part is holomorphically linearizable if its eigenvalue $\lambda\in\br$ satisfies the condition
\[
\sum_{n\geqslant 1}\frac{\log q_{n+1}}{q_n}<\infty,
\]
where $p_n/q_n$ is the $n^{\sl th}$ approximant of the continuous fraction expansion of $\lambda$.
\end{Thm}

In this article, all singularities are assumed to be linearizable. Let $(\bd^2,\cf,\{0\})$ be a holomorphic foliation on the unit bidisc $\bd^2$ defined by the linear vector field $Z=z\frac{\partial}{\partial z}+\lambda w\frac{\partial}{\partial w}$ with $\lambda\in\br^*$. One may assume $0<|\lambda|\leqslant 1$ after switching $z$ and $w$ if necessary. There are always two separatrices $\{z=0\}$ and $\{w=0\}$. Other leaves can be parametrized as
\[
L_\alpha:=\{(z,w)=\psi_\alpha(\zeta):=
(e^{i\,\zeta},\alpha\,e^{i\,\lambda\,\zeta})=(e^{-v+i\,u},\alpha\,e^{-\lambda\,v+i\,\lambda\,u})\}\cond{\alpha\neq 0},
\]
where $\zeta=u+iv\in\bc$. The map
\[
\aligned
\Psi:\bc\times\bc^*&\lrw\bc^2,\\
(\zeta,\alpha)&\longmapsto (e^{i\,\zeta},\alpha\,e^{i\,\lambda\,\zeta}),
\endaligned
\]
is locally biholomorphic. Here $\alpha$ is the coordinate on the transversal and $\zeta$ is the coordinate on leaves. It is not injective since $\Psi(\zeta+2\pi,\alpha)=\Psi(\zeta,\alpha\,e^{2\pi i\lambda})$.

Two numbers $\alpha$, $\beta\in\bc^*$ are {\sl equivalent} $\alpha\sim\beta$ if $\beta=e^{2k\pi i \lambda}\alpha$ for some $k\in\bz$. The following statements are equivalent:
\begin{itemize}
\item $\alpha\sim\beta$;
\item $L_\alpha=L_\beta$;
\item $\psi_\alpha=\psi_\beta\circ(\text{translation of }2k\pi)$ for some $k\in\bz$.
\end{itemize}

Let $\cc_\cf$ (resp. $\cc_\cf^{1,1}$) denote the space of functions \big(resp. forms of bidegree (1,1)\big) defined on leaves of the foliation which are compactly supported on $M\backslash E$, leafwise smooth and transversally continuous. A form $\iota\in\cc_\cf^{1,1}$ is said to be {\sl positive} if its restriction to every plaque is a positive (1,1)-form.

A {\sl directed harmonic current $T$ on $\cf$} is a continuous linear form on $\cc_\cf^{1,1}$ satisfying the following two conditions:
\begin{enumerate}
\item $i\partial\bar{\partial} T=0$ in the weak sense, i.e. $T(i\partial\bar{\partial}f)=0$ for all $f\in\cc_\cf$, where in the expression $i\partial\bar{\partial f}$, one only considers $\partial\bar{\partial}$ along the leaves;
\item $T$ is positive, that is $T(\iota)\geqslant 0$ for all positive forms $\iota\in\cc_\cf^{1,1}$.
\end{enumerate}

According to \cite{Berndtsson-Sibony-2002}, a directed harmonic current $T$ on a flow box $\bu\cong\bb\times\bt$ can be locally expressed as
\begin{equation}\label{local-expression}
T=\int_{\alpha\in\bt} h_\alpha [P_\alpha] d\mu(\alpha).
\end{equation}
The $h_\alpha$ are non-negative harmonic functions on the local leaves $P_\alpha$ and $\mu$ is a Borel measure on the transversal $\bt$. If $h_\alpha=0$ at some point on $P_\alpha$, then by the mean value theorem $h_\alpha\equiv 0$. For all such $\alpha\in\bt$, we replace $h_\alpha$ by the constant function $1$ and we set $d\mu(\alpha)=0$. Thus we get a new expression of $T$ where $h_\alpha>0$ for all $\alpha\in\bt$.

Such an expression is not unique since $T=\int_{\alpha\in\bt}\big(h_\alpha\,g(\alpha)\big)[P_\alpha]\big(\frac{1}{g(\alpha)}\,d\mu(\alpha)\big)$ for any bounded positive function $g:\bt\rw \br_{>0}$. The expression is unique after {\sl normalization}, which means that for each $\alpha\in\bt$ one fixes $h_\alpha(z_0,w_0)=1$ at some point $(z_0,w_0)\in P_\alpha$.

Each harmonic function $h_\alpha$ on the leaf $V_\alpha$ can be pulled back by the parametrization $\Psi$ as the harmonic function
\[
H_\alpha(u,v):=h_\alpha\big(e^{-v+iu},\alpha\,e^{-\lambda v+i\lambda u}\big).
\]
The domain of definition for $u$, $v$ will be precisely described later in this section.

In Section \ref{sect-intro} the notion of {\sl periodic current} was introduced. Here is an equivalent characterization.
\begin{Prop} A directed harmonic current $T$ is {\sl periodic} if and only if there exists some $k\in\bz_{>0}$ such that $H_\alpha(u+2 k\pi,v)=H_\alpha(u,v)$ for all $u,v$ and for $\mu$-almost all $\alpha$.
\end{Prop}
\proof By definition $T$ is invariant under $(z,w)\mapsto(z,e^{2k\pi i \lambda}w)$ for some $k\in\bz_{>0}$, which is equivalent to $H_\alpha(u+2 k\pi,v)=H_\alpha(u,v)$ for all $u,v$ and $\mu$-almost all $\alpha$.
\endproof

A current $T$ of the form \eqref{local-expression} is $dd^c$-closed on $\bd^2\backslash\{0\}$. But its trivial extension $\tilde{T}$ across the singularity $0$ is not necessarily $dd^c$-closed on $\bd^2$. It is true when $T$ is compactly supported, for example when $T$ is a localization of a current on a compact manifold, by the following argument of Dinh-Nguy\^en-Sibony in \cite{Dinh-Nguyen-Sibony-2012} Lemma 2.5.

Let $T$ be a directed harmonic current on $M\backslash E$, where $M$ is a compact complex manifold and $E$ is a finite set. The current $T$ can be extended by $0$ as a positive current $\tilde{T}$ on $M$. Next, apply

\begin{Thm}[Alessandrini-Bassanelli \cite{Alessandrini-Bassanelli-1993}] Let $\Omega$ be an open subset of $\bc^n$ and $Y$ an analytic subset of $\Omega$ of dimension $<p$. Suppose $T$ is a {\em negative} current of bidimension $(p,p)$ on $\Omega\backslash Y$ such that $dd^c T\geqslant 0$. Suppose $\tilde{T}$ is the trivial extension of $T$ across $Y$ by $0$. Then $dd^c\tilde{T}\geqslant 0$ on $\Omega$.
\end{Thm}

Here $-T$ is a negative current of bidimension $(1,1)$ on $M\backslash E$ with $dd^c (-T)\geqslant 0$ and $E$ has dimension $0$. So for the trivial extension $\tilde{T}$ on $M$ one has $dd^c(-\tilde{T})\geqslant 0$. Moreover $\tilde{T}$ is compactly supported since $M$ is compact. Thus
\[
\left<dd^c\tilde{T},1\right>=\left<\tilde{T},dd^c1\right>=0.
\]
Combining with $dd^c\tilde{T}\leqslant 0$ from the extension theorem, one concludes that $dd^c\tilde{T}=0$ on $M$. Thus locally near any singularity, the trivial extension $\tilde{T}$ is $dd^c$-closed.

Let $\beta:=idz\wedge d\bar{z}+idw\wedge d\bar{w}$ be the standard K\"ahler form. The {\sl mass} of $T$ on a domain $U\subset\bd^2$ is denoted by $|\!|T|\!|_U:=\int_U T\wedge \beta$. In this article, all currents are assumed to have finite mass on $\bd^2$.

\begin{Def} (See \cite[Subsection~2.4]{Nguyen-2020-09}){\bf .} Let $T$ be a directed harmonic current on $(\bd^2,\cf,\{0\})$. We define the {\sl Lelong number} by the limit
\[
\cl(T,0)=\limsup\limits_{r\rw0+}\frac{1}{\pi r^2}|\!|T|\!|_{r\bd^2}\in[0,+\infty].
\]
\end{Def}

The limit can be infinite when the trivial extension $\tilde{T}$ across the origin is not $dd^c$-closed \cite[Example~2.11]{Nguyen-2020-09}. When $\tilde{T}$ is $dd^c$-closed, the following theorem assures the finiteness.

\begin{Thm}[Skoda \cite{Skoda-1982}]\label{thm-skoda} Let $T$ be a positive $dd^c$-closed $(1,1)$-current in $\bd^2$. Then the function $r\mapsto \frac{1}{\pi r^2}|\!|T|\!|_{r\bd^2}$ is increasing with $r\in(0, 1]$.
\end{Thm}

In our case, the function
\[
r\mapsto \frac{1}{\pi r^2}|\!|\tilde{T}|\!|_{r\bd^2}=\frac{1}{\pi r^2}|\!|T|\!|_{r\bd^2},
\]
is increasing with $r\in(0, 1]$. In particular
\[
\cl(T,0)=\lim\limits_{r\rw0+}\frac{1}{\pi r^2}|\!|T|\!|_{r\bd^2}\in[0,\tfrac{1}{\pi}|\!|T|\!|_{\bd^2}].
\]

To calculate $|\!|T|\!|_{\bd^2}$ and $\cl(T,0)$ one shall study $\Psi^{-1}(r\,\bd^2)$ for $r\in(0,1]$. Define $P_\alpha:=L_\alpha\cap \bd^2$ and $P_\alpha^{(r)}:=L_\alpha\cap r\,\bd^2$. Define $\log^+(x):=\max\{0,\log(x)\}$ for $x>0$.

\begin{Lem} The range of $(u,v)$ for a point $(z,w)\in P_\alpha$ and $P_{\alpha}^{(r)}$ is either an upper-half plane or a horizontal strip. More precisely,
\begin{enumerate}
\item when $\lambda>0$,
\begin{align*}
(z,w)\in P_\alpha & \Longleftrightarrow  v>\frac{\log^+|\alpha|}{\lambda},\\
(z,w)\in P_\alpha^{(r)} & \Longleftrightarrow
\left\{
\begin{aligned}
&v>\frac{\log|\alpha|-\log r}{\lambda} &\cond{|\alpha|\geqslant r^{1-\lambda}},\\
&v>-\log r &\cond{|\alpha|<r^{1-\lambda}}.
\end{aligned}
\right.
\end{align*}

\item when $\lambda<0$, $P_\alpha=\emptyset$ for $|\alpha|\geqslant 1$, $P_\alpha^{(r)}=\emptyset$ for $|\alpha|\geqslant r^{1-\lambda}$ and for the other $\alpha$
\begin{align*}
(z,w)\in P_\alpha & \Longleftrightarrow  0<v<\frac{\log|\alpha|}{\lambda},\\
(z,w)\in P_\alpha^{(r)} & \Longleftrightarrow -\log r<v<\frac{\log|\alpha|-\log r}{\lambda}.
\end{align*}
\end{enumerate}
\end{Lem}
\proof
Recall that $(z,w)=(e^{-v+i\,u},\alpha\,e^{-\lambda\,v+i\,\lambda\,u})$ on $L_{\alpha}$. So for any $r\in(0,1]$, $(z,w)\in P_\alpha^{(r)}$ if and only if both $|z|=e^{-v}<r$ and $|w|=|\alpha|\,e^{-\lambda\,v}<r$.

When $\lambda>0$ one has $v>-\log r$ and $v>\frac{\log|\alpha|-\log r}{\lambda}$. In particular for $r=1$, one has $v>0$ and $v>\frac{\log|\alpha|}{\lambda}$.

When $\lambda<0$ one has $-\log r<v<\frac{\log|\alpha|-\log r}{\lambda}$. In particular for $r=1$, one has $0<v<\frac{\log|\alpha|}{\lambda}$. If there is no solution for $v$ then $P_{\alpha}^{(r)}=\emptyset$.\endproof

When $\lambda>0$, the range of $v$ is unbounded for each fixed $\alpha\in\bc^*$.

\medskip

\begin{minipage}{\linewidth}
      \centering
      \begin{minipage}{0.49\linewidth}
          \begin{figure}[H]
              \includegraphics[width=\linewidth]{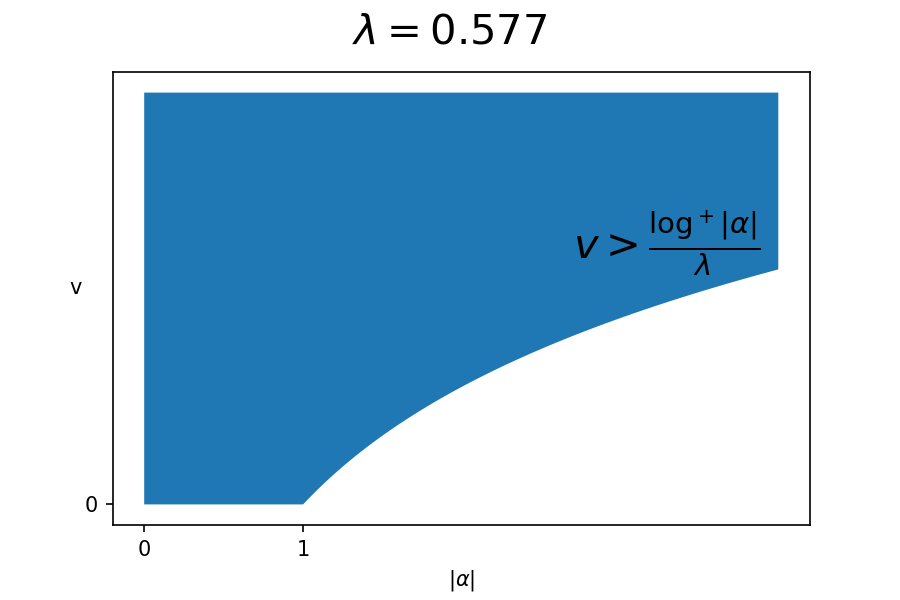}
              \caption{The region of $(|\alpha|,v)$ for $P_\alpha$}
          \end{figure}
      \end{minipage}
      \begin{minipage}{0.49\linewidth}
          \begin{figure}[H]
              \includegraphics[width=\linewidth]{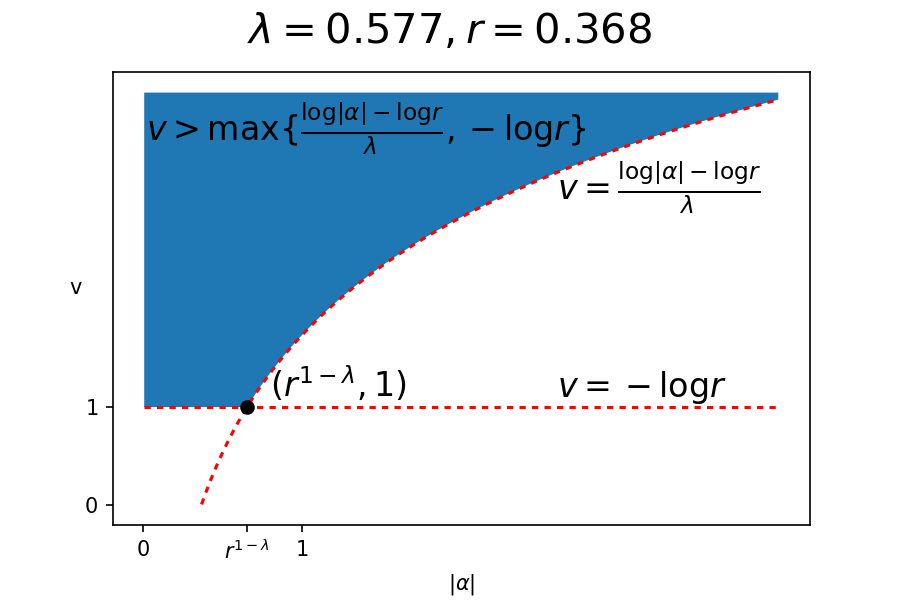}
              \caption{The region of $(|\alpha|,v)$ for $P_\alpha^{(r)}$}
          \end{figure}
      \end{minipage}
\end{minipage}

\medskip

When $\lambda<0$, the range of $v$ is bounded for each fixed $\alpha$.

\begin{minipage}{\linewidth}
      \centering
      \begin{minipage}{0.49\linewidth}
          \begin{figure}[H]
              \includegraphics[width=\linewidth]{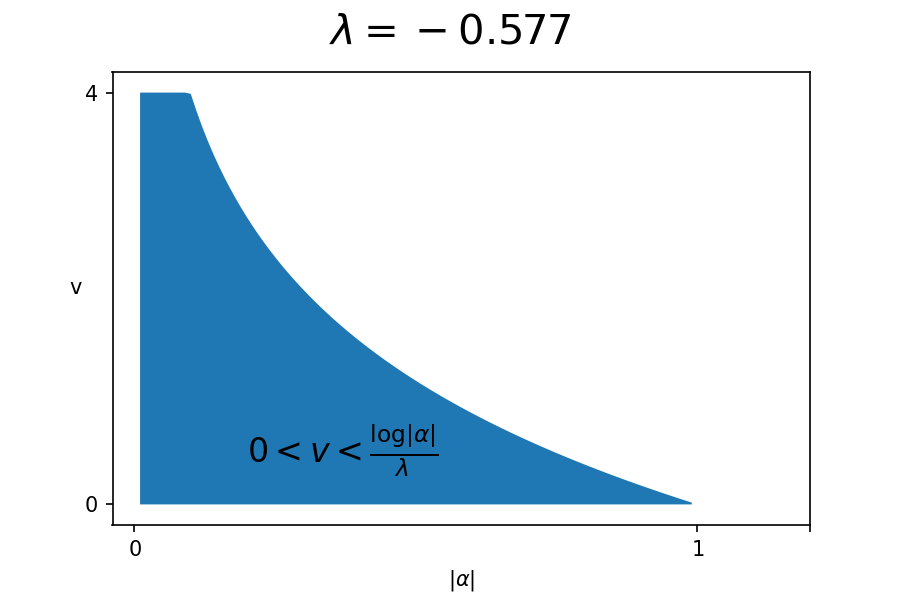}
              \caption{The region of $(|\alpha|,v)$ for $P_\alpha$}
          \end{figure}
      \end{minipage}
      \begin{minipage}{0.49\linewidth}
          \begin{figure}[H]
              \includegraphics[width=\linewidth]{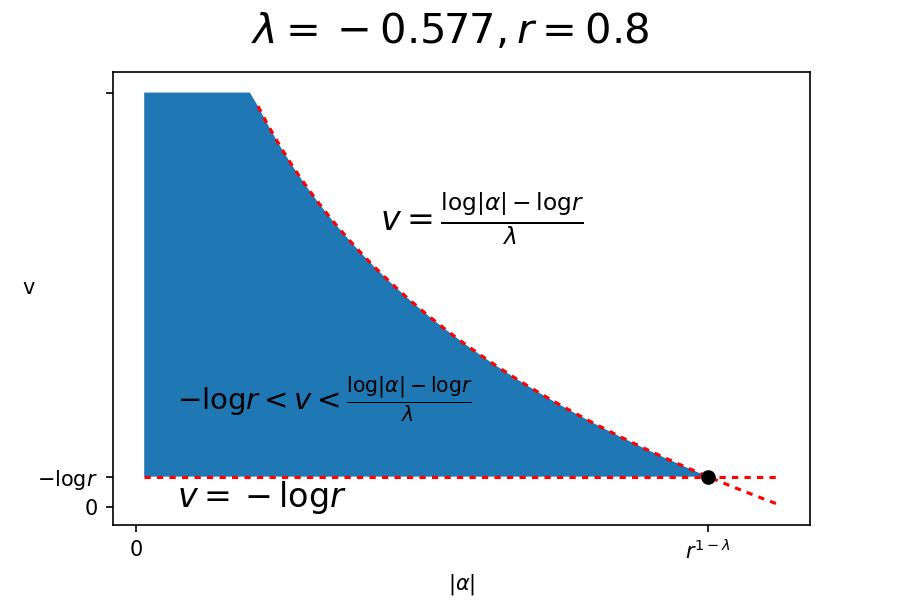}
              \caption{The region of $(|\alpha|,v)$ for $P_\alpha^{(r)}$}
          \end{figure}
      \end{minipage}
\end{minipage}

\medskip
In this article, the notations $\lesssim$ and $\gtrsim$ stand for inequalities up to a multiplicative non-zero constant depending only on $\lambda$. We write $\approx$ when both inequalities are satisfied.

\section{Geometry of leaves when $\lambda>0$}
For any $\alpha\in\bc^*$ fixed, the leaf $L_\alpha$ is contained in a real 3-dimensional Levi flat CR manifold $|w|=|\alpha|\,|z|^{\lambda}$, which can be viewed as a curve in $|z|=e^{-v}$, $|w|=|\alpha|\,e^{-\lambda v}$ coordinates. The norms $|z|$ and $|w|$ depends only on $v$. When $v\rw+\infty$, the point on the leaf tends to the singularity $(0,0)$ described by the following figures:

\begin{minipage}{\linewidth}
      \centering
      \begin{minipage}{0.49\linewidth}
          \begin{figure}[H]
              \includegraphics[width=\linewidth]{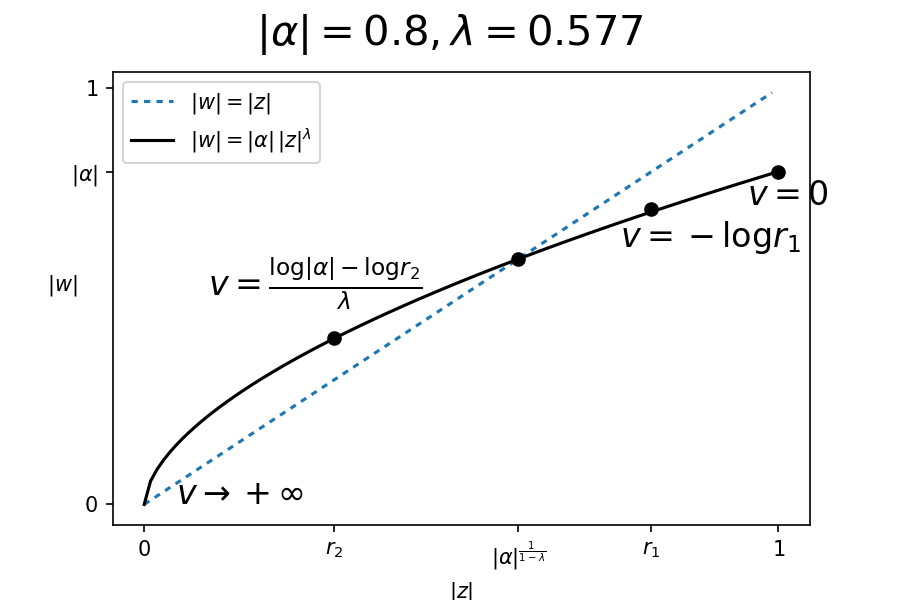}
              \caption{Case $|\alpha|<1$}
          \end{figure}
      \end{minipage}
      \begin{minipage}{0.49\linewidth}
          \begin{figure}[H]
              \includegraphics[width=\linewidth]{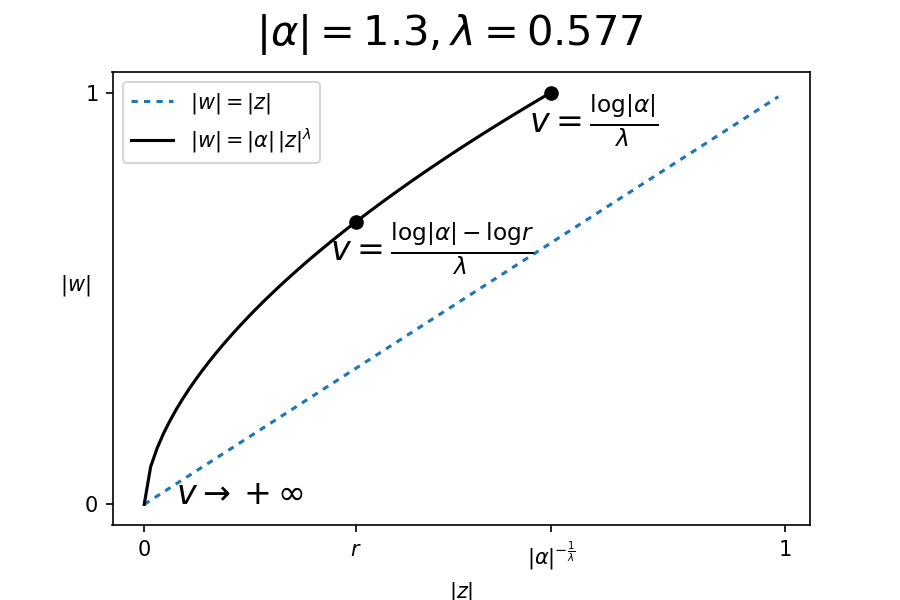}
              \caption{Case $|\alpha|\geqslant 1$}
          \end{figure}
      \end{minipage}
\end{minipage}

\medskip
If one fixes some $v=-\log r$, then $|z|=r$ and $|w|=|\alpha|\,r^\lambda$ is fixed. The set $\bt^2_r:=\{(z,w)\in\bd^2~:~|z|=r,|w|=|\alpha|\,r^\lambda\}$ is a torus and the intersection of the leaf $L_\alpha$ with this torus is a smooth curve $L_{\alpha,r}:=L_\alpha\cap\bt^2_r$.

When $\lambda\in\bq$, this curve $L_{\alpha,r}$ is closed.

\begin{figure}[H]
\begin{center}
   \includegraphics[width=0.45\linewidth]{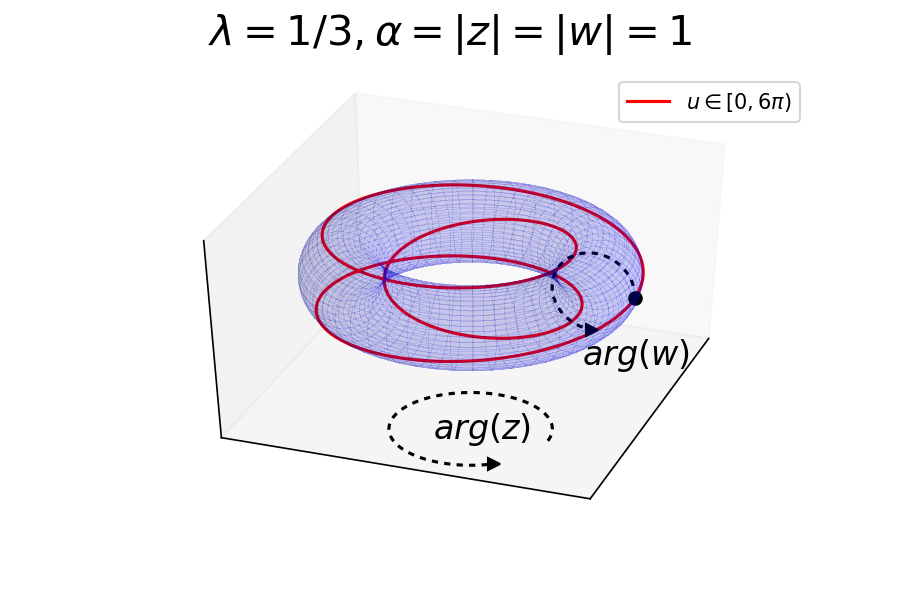}
   \caption{a closed curve on a torus}
   \end{center}
\end{figure}

When $\lambda\notin\bq$, this curve $L_{\alpha,r}$ is dense on the torus $\bt_r^2$.

\begin{minipage}{\linewidth}
      \centering
      \begin{minipage}{0.49\linewidth}
          \begin{figure}[H]
              \includegraphics[width=\linewidth]{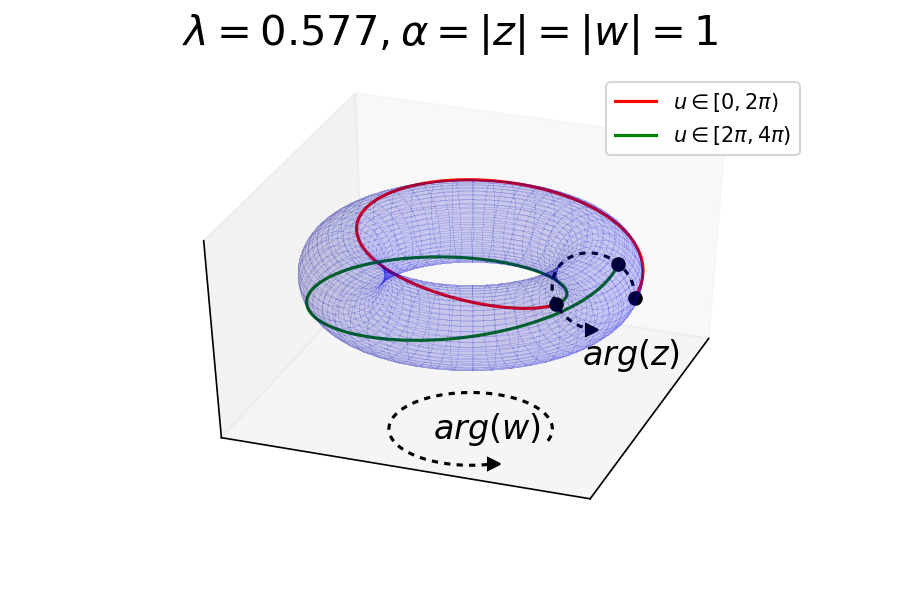}
              \caption{2 loops}\label{green}
          \end{figure}
      \end{minipage}
      \begin{minipage}{0.49\linewidth}
          \begin{figure}[H]
              \includegraphics[width=\linewidth]{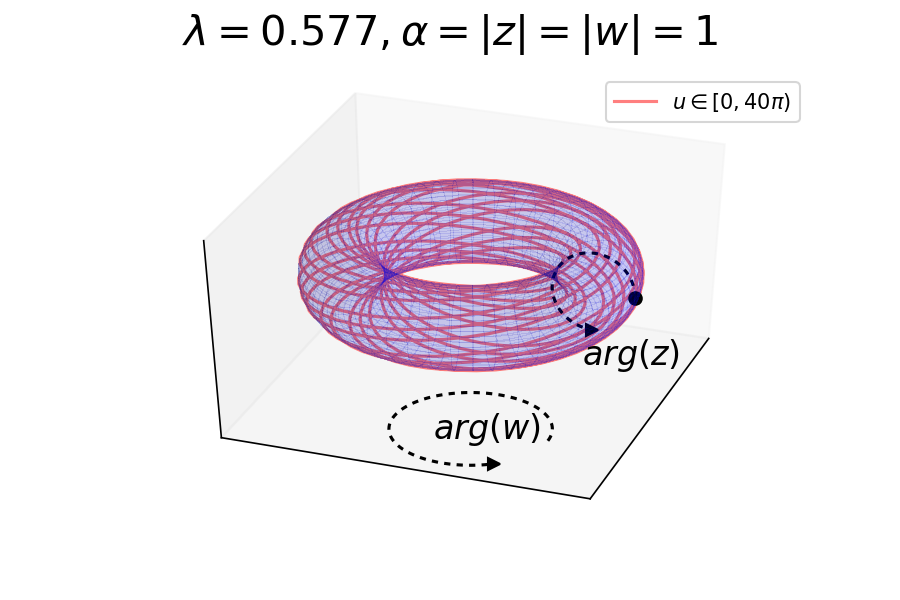}
              \caption{20 loops}
          \end{figure}
      \end{minipage}
\end{minipage}

\medskip
In this case the two curves $L_{\alpha,r}$ and $L_{\alpha\,e^{2\pi i \lambda},r}$ are two different parametrizations of the same image. The green curve in Figure \ref{green} is not only the image of $L_{\alpha,r}$ for $u\in[2\pi,4\pi)$ but also the image of $L_{\alpha\,e^{2\pi i \lambda},r}$ for $u\in[0,2\pi)$. This raises ambiguity while doing normalization on a leaf $L_\alpha$.

Such ambiguity can be resolved once one restricts everything to an open subset $U_\epsilon:=\{(z,w)\in\bd^2~|~{\rm arg}(z)\in(0,2\pi-\epsilon),z\neq 0,w\neq 0\}$ for some fixed $\epsilon\in[0,\pi)$. Any leaf $L_\alpha$ on $U_\epsilon$ decomposes into a disjoint union of infinitely many components:
\[
L_\alpha\cap U_\epsilon=\bigcup\limits_{k\in\bz}\big\{(e^{-v+iu},\alpha\,e^{2k\pi i\lambda}\,e^{-\lambda v+i\lambda u})~|~u\in(0,2\pi-\epsilon),v>\textstyle\frac{\log^+|\alpha|}{\lambda}\big\}.
\]
Such a parametrization is yet not unique. For example for any $k_0\in\bz$ one can do
\[
L_\alpha\cap U_\epsilon=\bigcup\limits_{k\in\bz}\big\{(e^{-v+iu},\alpha\,e^{2k\pi i\lambda}\,e^{-\lambda v+i\lambda u})~|~u\in(2k_0\pi,2k_0\pi+2\pi-\epsilon),v>\textstyle\frac{\log^+|\alpha|}{\lambda}\big\}.
\]
The parametrization is unique once one fixes $k_0$, for example $k_0=0$. But I leave a hint that all other choices of $k_0$ will be useful in a later proof.

Once one specifies the range of $u$ by fixing $k_0=0$, then the leaf on $U_\epsilon$ is uniquely parametrized.
\begin{figure}[H]
\begin{center}
   \includegraphics[width=0.45\linewidth]{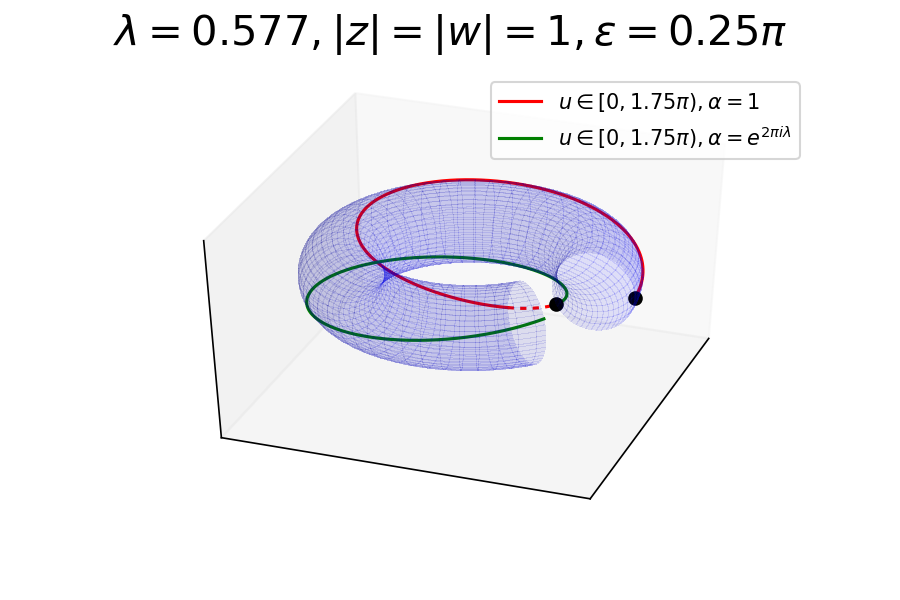}
   \caption{The curve $L_{1,1}$ in red and the curve $L_{e^{2\pi i \lambda},1}$ in green}
   \end{center}
\end{figure}
 
Later the limit case when $\epsilon=0$, $U:=\{(z,w)\in\bd^2~|~{\rm arg}(z)\in(0,2\pi),z\neq 0,w\neq 0\}=\{(z,w)\in\bd^2~|~z\notin\br_{\geqslant0},w\neq 0\}$ will be used.

\section{Rational case: $\lambda=\frac{a}{b}\in\bq$, $\lambda\in(0,1]$}
Say $\lambda=\frac{a}{b}$ where $a,b\in\bz_{\geqslant 1}$ are coprime. Then in $\bd^2$, for any $\alpha\in\bc^*$, the union $L_\alpha\cup\{0\}$ is the algebraic curve $\{w^b=\alpha^b\,z^a\}\cap\bd^2$. In other words, every leaf is a separatix. In this section it will be shown that any directed harmonic current $T$ has non-zero Lelong number.

The parametrization map $\psi_\alpha(\zeta):=(e^{i\zeta},\alpha\,e^{i\lambda\zeta})$ is now periodic: $\psi_\alpha(\zeta+2\pi b)=\psi_\alpha(\zeta)$. Let $T$ be a directed harmonic current. Then $T|_{P_\alpha}$ has the form $h_\alpha(z,w)[P_\alpha]$. Let 
\[
H_\alpha(u+iv):=h_\alpha\circ\psi_\alpha\big(u+iv+i\frac{\log^+|\alpha|}{\lambda}\big).
\]
It is a positive harmonic function for $\mu$-almost all $\alpha\in\bc^*$ defined in a neighborhood of the upper half plane $\bh:=\{(u+iv)\in\bc~|~v>0\}$. Moreover it is periodic: $H_\alpha(u+iv)=H_\alpha(u+2\pi b+iv)$. Periodic harmonic functions can be characterized by
\begin{Lem}\label{lem:periodic} Let $F(u,v)$ be a harmonic function in a neighborhood of $\bh$. If $F(u,v)=F(u+2\pi b,v)$ for all $(u,v)\in\bh$, then
\[
F(u,v)=\sum\limits_{k\in\bz,k\neq 0}\big(a_k\,e^{\frac{kv}{b}}\cos(\textstyle\frac{ku}{b})+b_k\,e^{\frac{kv}{b}}\sin(\frac{ku}{b})\big)+a_0+b_0\,v,
\]
for some $a_k$, $b_k\in\br$. Moreover, if $F|_\bh\geqslant 0$, then $a_0,b_0\geqslant 0$.
\end{Lem}
\proof
By periodicity
\[
F(u,v)=\sum\limits_{k=1}^{\infty}\big(A_k(v)\,\cos(\textstyle\frac{ku}{b})+B_k(v)\,\sin(\frac{ku}{b})\big)+A_0(v),
\]
for some functions $A_k(v)$, $B_k(v)$. They are smooth since $F$ is harmonic. Moreover
\[
0=\Delta F(u,v)=\sum\limits_{k=1}^{\infty}\Big(\big(A_k''(v)-(\textstyle\frac{k}{b})^2\,A_k(v)\big)\,\cos(\frac{ku}{b})+\big(B_k''(v)-(\frac{k}{b})^2\,B_k(v)\big)\,\sin(\frac{ku}{b})\Big)+A_0''(v).
\]
Thus
\[
A_k''(v)=(\textstyle\frac{k}{b})^2\,A_k(v), \ \ \ \ B_k''(v)=(\frac{k}{b})^2\,B_k(v), \ \ \ \ A_0''(v)=0.
\]
Hence
\[
A_k(v)=a_k\,e^{\frac{kv}{b}}+a_{-k}\,e^{-\frac{kv}{b}}, \ \ \ \ B_k(v)=b_k\,e^{\frac{kv}{b}}-b_{-k}\,e^{-\frac{kv}{b}}, \ \ \ \ A_0(v)=a_0+b_0\,v,
\]
for some $a_k$, $a_{-k}$, $b_k$, $b_{-k}\in\br$. One obtains the equality.

If $F|_\bh\geqslant 0$, then for any $v\geqslant 0$,
\[
\int_{u=0}^{2\pi b}F(u,v)du=2\pi b(a_0+b_0\,v)\geqslant0.
\]
Thus $a_0,b_0\geqslant 0.$
\endproof

For  $\alpha,\beta\in\bc^*$, the two maps $\psi_\alpha$ and $\psi_\beta$ parametrize the same leaf $L_\alpha=L_\beta$ if and only if $\beta=\alpha\,e^{2\pi i \frac{k}{b}}$ for some $k\in\bz$, i.e. $\alpha$ and $\beta$ differ from multiplying a $b^{\sl th}$-root of unity. Thus a transversal can be chosen as the section $\bt:=\{\alpha\in\bc^*~|~{\rm arg}(\alpha)\in[0,\frac{2\pi}{b})\}$. One fixes a normalization by setting $H_\alpha(0)=h_\alpha\circ \psi_\alpha(i\frac{\log^+|\alpha|}{\lambda})=1$.

The mass of the current $T$ is
\[
|\!\!|T|\!\!|_{\bd^2}=\int_{(z,w)\in\bd^2}T\wedge i\partial\bar{\partial}(|z|^2+|w|^2).
\]
In particular, one calculates the $(1,1)$-form $i\partial\bar{\partial}(|z|^2+|w|^2)$
on $L_\alpha$, where $z=e^{-v+iu},w=\alpha\,e^{-\lambda v+i\lambda u}$, using
\[
\aligned
dz&=ie^{-v+iu}du-e^{-v+iu}dv, & d\bar{z}&=-ie^{-v-iu}du-e^{-v-iu}dv,\\
dw&=i\alpha\,\lambda\,e^{-\lambda\,v+i\lambda\,u}du-\alpha\,\lambda\,e^{-\lambda\,v+i\lambda\,u}dv, & d\bar{w}&=-i\bar{\alpha}\,\lambda\,e^{-\lambda\,v-i\lambda\,u}du-\bar{\alpha}\,\lambda\,e^{-\lambda\,v-i\lambda\,u}dv,
\endaligned
\]
whence
\[
\aligned
i\partial\bar{\partial}(|z|^2+|w|^2)&=idz\wedge d\bar{z}+idw\wedge d\bar{w}\\
&=2\big(e^{-2v}+\lambda^2\,|\alpha|^2\,e^{-2\lambda\,v}\big)du\wedge dv.
\endaligned
\]
Thus
\[
\aligned
|\!|T|\!|_{\bd^2}&=\int_{\alpha\in\bt}h_\alpha(z,w)\int_{P_\alpha}i\partial\bar{\partial}(|z|^2+|w|^2)\,d\mu(\alpha)\\
&=\int_{\alpha\in\bt}\int_{u=0}^{2\pi b}\int_{v>0}H_\alpha(u+iv)\,2\big(e^{-2(v+\frac{\log^+|\alpha|}{\lambda})}+\lambda^2\,|\alpha|^2\,e^{-2\lambda\,(v+\frac{\log^+|\alpha|}{\lambda})}\big)\,du\wedge dv\,d\mu(\alpha)\\
&=\int_{\alpha\in\bt,|\alpha|<1}\int_{u=0}^{2\pi b}\int_{v>0}H_\alpha(u+iv)\,2\big(e^{-2v}+\lambda^2\,|\alpha|^2\,e^{-2\lambda\,v}\big)\,du\wedge dv\,d\mu(\alpha)\\
& \ \ \ +\int_{\alpha\in\bt,|\alpha|\geqslant1}\int_{u=0}^{2\pi b}\int_{v>0}H_\alpha(u+iv)\,2\big(|\alpha|^{-\frac{2}{\lambda}}\,e^{-2v}+\lambda^2\,e^{-2\lambda\,v}\big)du\wedge dv\,d\mu(\alpha).
\endaligned
\]

By Lemma \ref{lem:periodic},
\begin{equation}\label{eqn-periodic}
H_\alpha(u+iv)=\sum\limits_{k\in\bz,k\neq 0}\big(a_k(
\alpha)\,e^{\frac{kv}{b}}\cos(\textstyle\frac{ku}{b})+b_k(\alpha)\,e^{\frac{kv}{b}}\sin(\frac{ku}{b})\big)+a_0(\alpha)+b_0(\alpha)\,v,
\end{equation}
where $a_0(\alpha)$, $b_0(\alpha)$ are positive for $\mu$-almost all $\alpha$. Thus
\[
\aligned
|\!|T|\!|_{\bd^2}&=2\pi b\Big\{\int_{\alpha\in\bt,|\alpha|<1}\int_{v>0}\big(a_0(\alpha)+b_0(\alpha) v\big)\,2(e^{-2v}+\lambda^2\,|\alpha|^2\,e^{-2\lambda\,v})\,dv\,d\mu(\alpha)\\
& \ \ \ \ \ \ \ +\int_{\alpha\in\bt,|\alpha|\geqslant1}\int_{v>0}\big(a_0(\alpha)+b_0(\alpha) v\big)\,2(|\alpha|^{-\frac{2}{\lambda}}\,e^{-2v}+\lambda^2\,e^{-2\lambda\,v})\,dv\,d\mu(\alpha)\Big\}\\
&=2\pi b\,\Big\{\int_{|\alpha|<1}a_0(\alpha)\,(1+|\alpha|^2\,\lambda)d\mu(\alpha)+\int_{|\alpha|\geqslant1}a_0(\alpha)\,(|\alpha|^{-\frac{2}{\lambda}}+\lambda)d\mu(\alpha)\\
& \ \ \ \ \ \ +\int_{|\alpha|<1}b_0(\alpha)\,\big(\tfrac{1}{2}+\tfrac{1}{2}|\alpha|^2\big)d\mu(\alpha)+\int_{|\alpha|\geqslant1}b_0(\alpha)\,\big(\tfrac{1}{2}+\tfrac{1}{2}|\alpha|^{-\frac{2}{\lambda}}\big)d\mu(\alpha)\Big\}\\
&\approx \int_{\alpha\in\bc^*}a_0(\alpha)\,d\mu(\alpha)+\int_{\alpha\in\bc^*}b_0(\alpha)\,d\mu(\alpha).
\endaligned
\]

The Lelong number can now be calculated as follows
\[
\aligned
\cl(T,0)&=\lim\limits_{r\rw 0+}\frac{1}{r^2}|\!|T|\!|_{r\bd^2}\\
&=\lim\limits_{r\rw 0+}\frac{1}{r^2}2\pi b\,\Big\{\int_{\alpha\in\bt,|\alpha|<r^{1-\lambda}}\int_{v>-\log r}\big(a_0(\alpha)+b_0(\alpha)v\big)\,2\,(e^{-2v}+\lambda^2\,|\alpha|^2\,e^{-2\lambda\,v})\,dv\,d\mu(\alpha)\\
& \ \ \ \ \ \ \ \ \ \ \ \ \ \ \ \ +\,\int_{\alpha\in\bt,r^{1-\lambda}\leqslant|\alpha|<1}\int_{v>\frac{\log|\alpha|-\log r}{\lambda}}\big(a_0(\alpha)+b_0(\alpha)v\big)\,2\,(e^{-2v}+\lambda^2\,|\alpha|^2\,e^{-2\lambda\,v})\,dv\,d\mu(\alpha)\\
& \ \ \ \ \ \ \ \ \ \ \ \ \ \ \ \ +\,\int_{\alpha\in\bt,|\alpha|\geqslant1}\int_{v>\frac{-\log r}{\lambda}}\big(a_0(\alpha)+b_0(\alpha)v\big)\,2\,(|\alpha|^{-\frac{2}{\lambda}}e^{-2v}+\lambda^2\,e^{-2\lambda\,v})\,dv\,d\mu(\alpha)\Big\}\\
&=\lim\limits_{r\rw 0+}2\pi b\,\Big\{\int_{\alpha\in\bt,|\alpha|<r^{1-\lambda}}a_0(\alpha)\,(1+\lambda\,|\alpha|^2\,r^{2\lambda-2})\,d\mu(\alpha)\\
& \ \ \ \ \ \ \ \ \ \ \ \ \ \ +\int_{\alpha\in\bt,|\alpha|\geqslant r^{1-\lambda}}a_0(\alpha)\,(|\alpha|^{-\frac{2}{\lambda}}\,r^{\frac{2}{\lambda}-2}+\lambda)\,d\mu(\alpha)\\
& \ \ \ \ \ \ \ \ \ \ \ \ \ \ +\int_{\alpha\in\bt,|\alpha|< r^{1-\lambda}}b_0(\alpha)\,\big(\tfrac{1}{2}+\tfrac{1}{2}|\alpha|^2\,r^{2\lambda-2}-\log r-\lambda\,|\alpha|^2\,r^{2\lambda-2}\,\log r\big)\,d\mu(\alpha)\\
& \ \ \ \ \ \ \ \ \ \ \ \ \ \ +\int_{\alpha\in\bt,r^{1-\lambda}\leqslant|\alpha|<1}b_0(\alpha)\,\big(\tfrac{1}{2}+\tfrac{1}{2}|\alpha|^{-\frac{2}{\lambda}}\,r^{\frac{2}{\lambda}-2}-\log r-|\alpha|^{-\frac{2}{\lambda}}\,\lambda^{-1}\,r^{2\lambda-2}\log r\\
& \ \ \ \ \ \ \ \ \ \ \ \ \ \ \ \ \ \ \ \ \ \ \ \ \ \ \ \ \ \ \ \ \ \ \ \ \ \ \ \ \ \ \ \ +\log |\alpha|+\lambda^{-1}\,|\alpha|^{-\frac{2}{\lambda}}\,\log|\alpha|\,r^{2\lambda-2}\big)\,d\mu(\alpha)\\
& \ \ \ \ \ \ \ \ \ \ \ \ \ \ +\int_{\alpha\in\bt,|\alpha|\geqslant 1}b_0(\alpha)\,\big(\tfrac{1}{2}+\tfrac{1}{2}|\alpha|^{-\frac{2}{\lambda}}\,r^{\frac{2}{\lambda}-2}-\log r-\lambda^{-1}\,|\alpha|^{-\frac{2}{\lambda}}\,r^{2\lambda-2}\log r\big)\,d\mu(\alpha)
\Big\}.
\endaligned
\]

First one analyzes the $a_0(\alpha)$ part. When $|\alpha|<r^{1-\lambda}$,
\begin{align}\label{ineq1}
1<1+\lambda\,|\alpha|^2\,r^{2\lambda-2}<1+\lambda\,r^{2-2\lambda}\,r^{2\lambda-2}=1+\lambda,
\end{align}
is uniformly bounded with respect to $\alpha$ and $r$. When $|\alpha|\geqslant r^{1-\lambda}$
\begin{align}\label{ineq2}
\lambda<|\alpha|^{-\frac{2}{\lambda}}\,r^{\frac{2}{\lambda}-2}+\lambda<1+\lambda,
\end{align}
is also uniformly bounded with respect to $\alpha$ and $r$. Thus
\[
\cl(T,0)\approx \underbrace{\int_{\alpha\in\bt}a_0(\alpha)\,d\mu(\alpha)}_{\text{linear part}}+\underbrace{\lim\limits_{r\rw 0+}\big(\text{$b_0(\alpha)$ part}\big)}_{\text{with $v$ part}}.
\]

Next one analyses the $b_0(\alpha)$ part.
\begin{Lem}\label{lem:no-C-alpha} The Lelong number of $T$ at $0$ is finite only if $b_0(\alpha)=0$ for $\mu$-almost all $\alpha\in\bt$.
\end{Lem}
\proof Suppose not, i.e. $\int_{\alpha\in\bt}b_0(\alpha)\,d\mu(\alpha)=B_0>0$. Then
\[
\aligned
\cl(T,0)&\geqslant\lim\limits_{r\rw 0+}2\pi b\Big\{\int_{\alpha\in\bt,|\alpha|< r^{1-\lambda}}b_0(\alpha)\,\big(-\log r\big)\,d\mu(\alpha)+\int_{\alpha\in\bt,|\alpha|\geqslant r^{1-\lambda}}b_0(\alpha)\,\big(-\log r\big)\,d\mu(\alpha)\Big\}\\
&=2\pi b\,B_0\,\lim\limits_{r\rw 0+}(-\log r)=+\infty,
\endaligned
\]
would contradict the finiteness of the Lelong number stated in Theorem \ref{thm-skoda}.
\endproof

Thus one may assume $b_0(\alpha)=0$ for $\mu$-almost all $\alpha\in\bt$. Then the Lelong number
\[
\cl(T,0)\approx\int_{\alpha\in\bt}a_0(\alpha)\,d\mu(\alpha)\approx|\!|T|\!|_{\bd^2},
\]
is strictly positive.

\section{Irrational case $\lambda\notin \bq$, $\lambda\in (0,1)$}
Now $\{z=0\}$ and $\{w=0\}$ are the only 2 separatrices in $\bd^2$. For each fixed $\alpha\in\bc^*$, the map $\psi_\alpha(\zeta)=(e^{i\,\zeta},\alpha\,e^{i\,\lambda\,\zeta})$ is injective since $\lambda\notin\bq$.

Any pair of equivalent numbers $\alpha\sim\beta$, $\beta=\alpha\,e^{2k\pi i \lambda}$, may provide us with two different normalizations $H_{\alpha}$ and $H_{\beta}$ on the same leaf $L_{\alpha}=L_{\beta}$. A major task is to find formulas for the mass and the Lelong number independent by the choice of normalization.

Let $T$ be a harmonic current directed by $\cf$. Then $T|_{P_\alpha}$ has the form $h_\alpha(z,w)[P_\alpha]$. One may assume that $h_\alpha$ is nowhere 0 for every $\alpha$. Let
\[
H_\alpha(u+iv):=h_\alpha\circ \psi_\alpha\big(u+iv+i\frac{\log^+|\alpha|}{\lambda}\big).
\]
 It is a positive harmonic function for $\mu$-almost all $\alpha\in\bc^*$ defined in a neighborhood of the upper half plane $\bh=\{(u+iv)\in\bc~|~v>0\}$, determined by the Poisson integral formula
\[
H_\alpha(u+iv)=\frac{1}{\pi}\,\int_{y\in\br}H_\alpha(y)\,\frac{v}
{v^2+(y-u)^2}\,dy+C_\alpha\,v.
\]
One can normalize $H_\alpha$ by setting $H_\alpha(0)=1$. But by doing so one may normalize data over the same leaf for multiple times. The ambiguity is described by

\begin{Prop}\label{prop:ab}
If $\beta=\alpha\,e^{2k\,\pi\,i\,\lambda}$ for some $k\in\bz$, then the two normalized positive harmonic functions $H_\alpha$ and $H_\beta$ satisfy
\[
H_\alpha(u+iv)
=
H_\alpha(2k\pi)\,H_\beta(u-2k\pi+iv).
\]
In other words, they differ by a translation and a multiplication by a non-zero constant.
\end{Prop}

\proof
When $|\alpha|<1$, by definition
\[
H_\alpha(u+iv)=h_\alpha(e^{-v+iu},\alpha\,e^{-\lambda\,v+i\,\lambda\,u}), \ \ \ \ H_\alpha(0)=h_\alpha(1,\alpha).
\]
Thus the normalized harmonic function is
\[
H_\alpha(u+iv)=\frac{h_\alpha(e^{-v+iu},\alpha\,e^{-\lambda\,v+i\,\lambda\,u})}{h_\alpha(1,\alpha)},
\]
and for the same reason
\[H_\beta(u+iv)=\frac{h_\beta(e^{-v+iu},\beta\,e^{-\lambda\,v+i\,\lambda\,u})}{h_\beta(1,\beta)}.
\]

The two functions $h_\alpha$ and $h_\beta$ are the positive harmonic coefficient of $T$ on the same leaf $L_\alpha=L_\beta$, hence they differ up to multiplication by a positive constant $C>0$,
\[
\aligned
h_\alpha(e^{-v+iu},\alpha\,e^{-\lambda\,v+i\,\lambda\,u})&=C\cdot h_\beta(e^{-v+iu},\alpha\,e^{-\lambda\,v+i\,\lambda\,u})\\
&=C\cdot h_\beta(e^{-v+iu},\beta\,e^{-2k\,\pi\,i\,\lambda}\,e^{-\lambda\,v+i\,\lambda\,u})\\
&=C\cdot h_\beta(e^{-v+i(u-2k\,\pi)},\beta\,e^{-\lambda\,v+i\,\lambda\,(u-2k\,\pi)}).
\endaligned
\]
Thus
\[
\aligned
H_\alpha(u+iv)&=\frac{h_\alpha(e^{-v+iu},\alpha\,e^{-\lambda\,v+i\,\lambda\,u})}{h_\alpha(1,\alpha)}=\frac{C\cdot h_\beta(e^{-v+i(u-2k\,\pi)},\beta\,e^{-\lambda\,v+i\,\lambda\,(u-2k\,\pi)})}{C\cdot h_\beta(1,\alpha)}\\
&=\frac{h_\beta(e^{-v+i(u-2k\,\pi)},\beta\,e^{-\lambda\,v+i\,\lambda\,(u-2k\,\pi)})}{h_\beta(1,\beta)}\cdot \frac{h_\beta(1,\beta)}{h_\beta(1,\alpha)}\\
&=H_\beta(u-2k\,\pi+iv)\cdot \frac{h_\beta(1,\beta)}{h_\beta(1,\alpha)}.
\endaligned
\]
When $u=2k\,\pi$ and $v=0$ one has $H_\alpha(2k\,\pi)=\frac{h_\beta(1,\beta)}{h_\beta(1,\alpha)}$. Thus one gets the equality. The proof for the case $|\alpha|>1$ is similar.
\endproof

Take the open subset $U:=\{(z,w)\in\bd^2~|~z\notin\br_{\geqslant0},w\neq 0\}$.

\begin{minipage}{\linewidth}
      \centering
      \begin{minipage}{0.49\linewidth}
          \begin{figure}[H]
              \includegraphics[width=\linewidth]{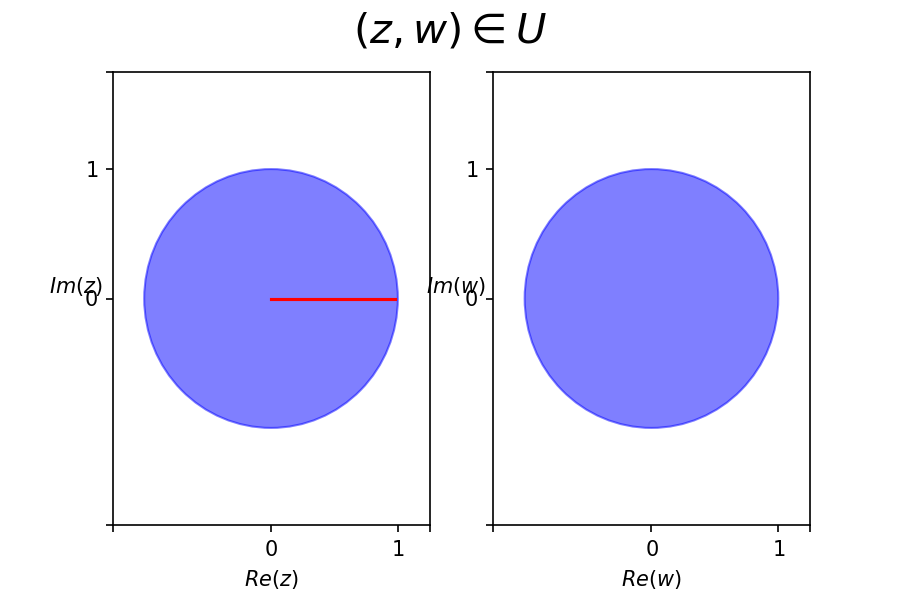}
              \caption{Domain $U$ in coordinates $(z,w)$}
          \end{figure}
      \end{minipage}
      \begin{minipage}{0.49\linewidth}
          \begin{figure}[H]
              \includegraphics[width=\linewidth]{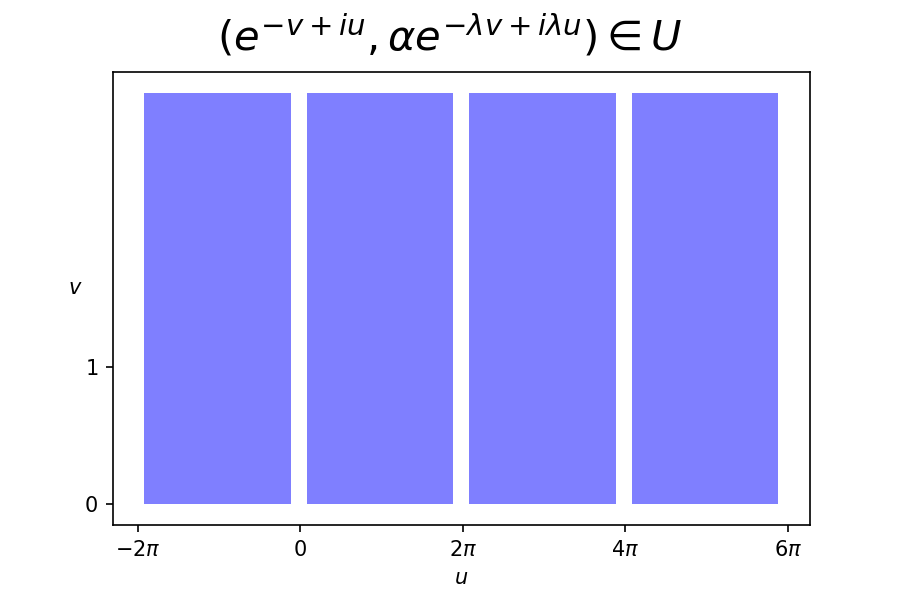}
              \caption{Domain $U$ in coordinates $(u,v)$}
          \end{figure}
      \end{minipage}
\end{minipage}

\medskip

Any leaf $L_\alpha$ in $U$ is a disjoint union of infinitely many components:
\[
L_\alpha\cap U=\bigcup\limits_{k\in\bz}\tilde{L}_{\alpha\,e^{2k\pi i\lambda}}:=\bigcup\limits_{k\in\bz}\big\{(e^{-v+iu},\alpha\,e^{2k\pi i\lambda}\,e^{-\lambda v+i\lambda u})~|~u\in(0,2\pi),v>\textstyle\frac{\log^+|\alpha|}{\lambda}\big\}.
\]

Normalizing $H_{\alpha\,e^{2k\pi i\lambda}}$ on $\tilde{L}_{\alpha\,e^{2k\pi i\lambda}}$ avoids ambiguity. Thus
\[
|\!|T|\!|_{U}=\int_{\alpha\in\bc^*}\int_{v>0}\int_{u=0}^{2\pi}H_{\alpha}(u+iv)|\!|\psi_\alpha'|\!|^2du\,dv\,d\mu(\alpha)
\]
for some positive measure $\mu$ on $\bc^*$. Here, let us write $|\!|\psi_\alpha'|\!|^2$ for the jacobian calculated before,
\begin{align}\label{jacobian}
|\!|\psi_\alpha'|\!|^2=\left\{
\begin{aligned}
&2(e^{-2v}+\lambda^2\,|\alpha|^2\,e^{-2\lambda\,v}), &\cond{|\alpha|< 1}\\
&2(|\alpha|^{-\frac{2}{\lambda}}e^{-2v}+\lambda^2\,e^{-2\lambda\,v}). &\cond{|\alpha|\geqslant 1}
\end{aligned}
\right.
\end{align}

Since $H$ is harmonic in a neighborhood of $\bh$, it is continuous in $\bh$. So
\[
\aligned
|\!|T|\!|_{U}&=\lim\limits_{\epsilon\rw 0+}\int_{\alpha\in\bc^*}\int_{v>0}\int_{u=0}^{2\pi+\epsilon}H_{\alpha}(u+iv)|\!|\psi_\alpha'|\!|^2du\,dv\,d\mu(\alpha)\\
&=\lim\limits_{\epsilon\rw0+}|\!|T|\!|_{\bigcup\limits_{k\in\bz}\tilde{L}_{\alpha\,e^{2k \pi i\lambda}}}\\
&=|\!|T|\!|_{\bd^2}.
\endaligned
\]
Thus we can express the mass by a formula independent of the choice of normalization
\[
|\!|T|\!|_{\bd^2}=\int_{\alpha\in\bc^*}\int_{v>0}\int_{u=0}^{2\pi}H_{\alpha}(u+iv)|\!|\psi_\alpha'|\!|^2du\,dv\,d\mu(\alpha).
\]

\begin{Lem}\label{Lem:Mass}
For each $k_0\in\bz$ fixed,
\[
|\!|T|\!|_{\bd^2}=\int_{\alpha\in\bc^*}\int_{v>0}\int_{u=2k_0\pi}^{2k_0\pi+2\pi}H_{\alpha}(u+iv)|\!|\psi_\alpha'|\!|^2du\,dv\,d\mu(\alpha).
\]
\end{Lem}
\proof
The disjoint union $L_\alpha\cap U=\bigcup\limits_{k\in\bz}\tilde{L}_{\alpha\,e^{2k\pi i\lambda}}$ can be parametrized in many other ways. For instance
\[
L_\alpha\cap U=\bigcup\limits_{k\in\bz}\big\{(e^{-v+iu},\alpha\,e^{2k\pi i\lambda}\,e^{-\lambda v+i\lambda u})~|~u\in(2k_0\pi,2k_0\pi+2\pi),v>\textstyle\frac{\log^+|\alpha|}{\lambda}\big\}.
\]
By the same argument as above one concludes.
\endproof

\subsection{Periodic currents, still a Fourrier series}
Periodic currents behave similarly as currents in the rational case $\lambda\in\bq$. Suppose $H_\alpha$ is periodic, i.e. there is some $b\in\bz_{\geqslant 1}$ such that $H_\alpha(u+iv)=H_\alpha(u+2\pi b+iv)$ for any $u+iv\in\bh$. Periodic harmonic functions are characterized as \eqref{eqn-periodic} of Lemma~\ref{lem:periodic}.

According to Lemma~\ref{Lem:Mass}, the mass is
\[
\aligned
|\!|T|\!|_{\bd^2}&=\int_{\alpha\in\bc^*}\int_{v>0}\int_{u=2k_0\pi}^{2k_0\pi+2\pi}H_\alpha(u+iv)|\!|\psi_\alpha'|\!|^2\,du\wedge dv\,d\mu(\alpha),
\endaligned
\]
for any $k_0\in\bz$, in particular, for $k_0=0,1,\dots,b-1$. Thus, we may calculate
\[
\aligned
b|\!|T|\!|_{\bd^2}&=\int_{\alpha\in\bc^*}\int_{v>0}\int_{u=0}^{2\pi b}H_\alpha(u+iv)|\!|\psi_\alpha'|\!|^2\,du\wedge dv\,d\mu(\alpha),\\
|\!|T|\!|_{\bd^2}&=\frac{1}{b}\int_{\alpha\in\bc^*}\int_{v>0}\int_{u=0}^{2\pi b}H_\alpha(u+iv)|\!|\psi_\alpha'|\!|^2\,du\wedge dv\,d\mu(\alpha),\\
&=\frac{1}{b}\Big\{\int_{|\alpha|<1}\int_{v>0}\int_{u=0}^{2\pi b}H_\alpha(u+iv)\,2(e^{-2v}+\lambda^2\,|\alpha|^2\,e^{-2\lambda\,v})\,du\wedge dv\,d\mu(\alpha)\\
& \ \ \ \ \ \ \ +\int_{|\alpha|\geqslant1}\int_{v>0}\int_{u=0}^{2\pi b}H_\alpha(u+iv)\,2(|\alpha|^{-\frac{2}{\lambda}}\,e^{-2v}+\lambda^2\,e^{-2\lambda\,v})\,du\wedge dv\,d\mu(\alpha)\Big\},\\
&=\frac{2\pi b}{b}\Big\{\int_{|\alpha|<1}\int_{v>0}\big(a_0(\alpha)+b_0(\alpha) v\big)\,2(e^{-2v}+\lambda^2\,|\alpha|^2\,e^{-2\lambda\,v})\,dv\,d\mu(\alpha)\\
& \ \ \ \ \ \ \ +\int_{|\alpha|\geqslant1}\int_{v>0}\big(a_0(\alpha)+b_0(\alpha) v\big)\,2(|\alpha|^{-\frac{2}{\lambda}}\,e^{-2v}+\lambda^2\,e^{-2\lambda\,v})\,dv\,d\mu(\alpha)\Big\},\\
&=2\pi \,\Big\{\int_{|\alpha|<1}a_0(\alpha)\,(1+|\alpha|^2\,\lambda)d\mu(\alpha)+\int_{|\alpha|\geqslant1}a_0(\alpha)\,(|\alpha|^{-\frac{2}{\lambda}}+\lambda)d\mu(\alpha)\\
& \ \ \ \ \ \ +\int_{|\alpha|<1}b_0(\alpha)\,\big(t\frac{1}{2}+\tfrac{1}{2}|\alpha|^2\big)d\mu(\alpha)+\int_{|\alpha|<1}b_0(\alpha)\,\big(\tfrac{1}{2}+\tfrac{1}{2}|\alpha|^{-\frac{2}{\lambda}}\big)d\mu(\alpha)\Big\}\\
&\approx \int_{\alpha\in\bc^*}a_0(\alpha)\,d\mu(\alpha)+\int_{\alpha\in\bc^*}b_0(\alpha)\,d\mu(\alpha),
\endaligned
\]
which is the same expression as in the case $\lambda\in\bq_{>0}$.

Next, the Lelong number is calculated as
\[
\aligned
\cl(T,0)&=\lim\limits_{r\rw 0+}\frac{1}{r^2}|\!|T|\!|_{r\bd^2}\\
&=\lim\limits_{r\rw 0+}\frac{1}{r^2}2\pi\,\Big\{\int_{|\alpha|<r^{1-\lambda}}\int_{v>-\log r}\big(a_0(\alpha)+b_0(\alpha)v\big)\,2\,(e^{-2v}+\lambda^2\,|\alpha|^2\,e^{-2\lambda\,v})\,dv\,d\mu(\alpha)\\
& \ \ \ \ \ \ \ \ \ \ \ \ \ \ \ \ +\,\int_{r^{1-\lambda}\leqslant |\alpha|<1}\int_{v>\frac{\log|\alpha|-\log r}{\lambda}}\big(a_0(\alpha)+b_0(\alpha)v\big)\,2\,(e^{-2v}+\lambda^2\,|\alpha|^2\,e^{-2\lambda\,v})\,dv\,d\mu(\alpha)\Big\}\\
& \ \ \ \ \ \ \ \ \ \ \ \ \ \ \ \ +\,\int_{|\alpha|\geqslant 1}\int_{v>\frac{-\log r}{\lambda}}\big(a_0(\alpha)+b_0(\alpha)v\big)\,2\,(|\alpha|^{-\frac{2}{\lambda}}e^{-2v}+\lambda^2\,e^{-2\lambda\,v})\,dv\,d\mu(\alpha)\Big\}\\
&=\lim\limits_{r\rw 0+}2\pi\,\Big\{\int_{|\alpha|<r^{1-\lambda}}a_0(\alpha)\,(1+\lambda\,|\alpha|^2\,r^{2\lambda-2})\,d\mu(\alpha)\\
& \ \ \ \ \ \ \ \ \ \ \ \ \  +\int_{|\alpha|\geqslant r^{1-\lambda}}a_0(\alpha)\,(|\alpha|^{-\frac{2}{\lambda}}\,r^{\frac{2}{\lambda}-2}+\lambda)\,d\mu(\alpha)\\
& \ \ \ \ \ \ \ \ \ \ \ \ \  +\int_{|\alpha|< r^{1-\lambda}}b_0(\alpha)\,\big(\tfrac{1}{2}+\tfrac{1}{2}|\alpha|^2\,r^{2\lambda-2}-\log r-\lambda\,|\alpha|^2\,r^{2\lambda-2}\,\log r\big)\,d\mu(\alpha)\\
& \ \ \ \ \ \ \ \ \ \ \ \ \  +\int_{r^{1-\lambda}\leqslant |\alpha|<1}b_0(\alpha)\,\big(\tfrac{1}{2}+\tfrac{1}{2}|\alpha|^{-\frac{2}{\lambda}}\,r^{\frac{2}{\lambda}-2}-\log r-\lambda^{-1}\,|\alpha|^{-\frac{2}{\lambda}}\,r^{2\lambda-2}\log r\\
& \ \ \ \ \ \ \ \ \ \ \ \ \ \ \ \ \ \ \ \ \ \ \ \ \ \ \ \ \ \ \ \ \ \ \ \ \ \ \ +\log |\alpha|+\lambda^{-1}\,|\alpha|^{-\frac{2}{\lambda}}\,\log|\alpha|\,r^{2\lambda-2}\big)\,d\mu(\alpha)\\
& \ \ \ \ \ \ \ \ \ \ \ \ \  +\int_{|\alpha|\geqslant 1}b_0(\alpha)\,\big(\tfrac{1}{2}+\tfrac{1}{2}|\alpha|^{-\frac{2}{\lambda}}\,r^{\frac{2}{\lambda}-2}-\log r-\lambda^{-1}\,|\alpha|^{-\frac{2}{\lambda}}\,r^{2\lambda-2}\log r\big)\,d\mu(\alpha)\Big\}.
\endaligned
\]
exactly the same expression as in the case $\lambda=1$. Using the same argument as in Lemma \ref{lem:no-C-alpha} one may assume that $b_0(\alpha)=0$ for $\mu$-almost all $\alpha\in\bc^*$. One concludes that
\[
\cl(T,0)\approx\int_{\alpha\in\bc^*}a_0(\alpha)\,d\mu(\alpha)\approx|\!|T|\!|_{\bd^2}.
\]
The Lelong number is strictly positive, the same as in the case $\lambda\in\bq\cup(0,1)$.

\subsection{Non-periodic current}
For periodic currents, one takes an average among $b$ expressions \eqref{lem:periodic} in the previous section. For non-periodic currents, there is no canonical way of normalization. The key technique is to calculate expressions \eqref{lem:periodic} for all $k_0\in\bz$.

The Lelong number is expressed by
\[
\aligned
\cl(T,0)&=\lim\limits_{r\rw 0+}\frac{1}{r^2}\Big\{\int_{|\alpha|<r^{1-\lambda}}\int_{v>-\log r}\int_{u=0}^{2\pi}H_{\alpha}(u+iv)|\!|\psi_\alpha'|\!|^2du\,dv\,d\mu(\alpha)\\
& \ \ \ \ \ \ \ \ \ \ \ \ \ +\int_{r^{1-\lambda}\leqslant|\alpha|<1}\int_{v>\frac{\log|\alpha|-\log r}{\lambda}}\int_{u=0}^{2\pi}H_{\alpha}(u+iv)|\!|\psi_\alpha'|\!|^2du\,dv\,d\mu(\alpha)\\
& \ \ \ \ \ \ \ \ \ \ \ \ \ +\int_{|\alpha|\geqslant1|}\int_{v>\frac{-\log r}{\lambda}}\int_{u=0}^{2\pi}H_{\alpha}(u+iv)|\!|\psi_\alpha'|\!|^2du\,dv\,d\mu(\alpha)\Big\}
\endaligned
\]

Recall the Poisson integral formula after multiplying a nonzero constant
\[
H_\alpha(u+iv)=\frac{1}{\pi}\,\int_{y\in\br}H_\alpha(y)\,\frac{v}
{v^2+(y-u)^2}\,dy+{C}_\alpha\,v.
\]
Using the same argument as in Lemma \ref{lem:no-C-alpha}, one may assume ${C}_\alpha=0$ for all $\alpha\in\bc^*$.

\begin{Lem}\label{approx-poission-lambda>0}
For any $v\geqslant \frac{1}{\lambda}>1$ and for any $u\in\br$,
\[
\aligned
\frac{\frac{\partial}{\partial v}\Big(-\frac{1}{2}\frac{v}{v^2+(u-y)^2}e^{-2v}\Big)}{\frac{v}{v^2+(u-y)^2}e^{-2v}}\in\Big(\frac{1}{2},2\Big), \\
\frac{\frac{\partial}{\partial v}\Big(-\frac{1}{2\lambda}\frac{v}{v^2+(u-y)^2}e^{-2\lambda v}\Big)}{\frac{v}{v^2+(u-y)^2}e^{-2\lambda v}}\in\Big(\frac{1}{2},2\Big).
\endaligned
\]
\end{Lem}
\proof
This can be calculated directly,
\[
\aligned
\frac{\partial}{\partial v}\Big(-\frac{1}{2}\frac{v}{v^2+(u-y)^2}e^{-2v}\Big)&=\Big(\frac{v}{v^2+(u-y)^2}+\big(-\frac{1}{2}\big)\frac{1}{v^2+(u-y)^2}+\big(-\frac{1}{2}\big)\frac{v(-2v)}{\big(v^2+(u-y)^2\big)^2}\Big)\,e^{-2v}\\
\frac{\frac{\partial}{\partial v}\Big(-\frac{1}{2}\frac{v}{v^2+(u-y)^2}e^{-2v}\Big)}{\frac{v}{v^2+(u-y)^2}e^{-2v}}&=1+\big(-\frac{1}{2}\frac{1}{v}\big)+\frac{v}{v^2+(u-y)^2}\\
&\in\Big(1-\frac{1}{2v},1+\frac{1}{v}\Big)\subseteq\Big(\frac{1}{2},2\Big),\cond{v>1},\\
\frac{\partial}{\partial v}\Big(-\frac{1}{2\lambda}\frac{v}{v^2+(u-y)^2}e^{-2\lambda v}\Big)&=\Big(\frac{v}{v^2+(u-y)^2}+\big(-\frac{1}{2\lambda}\big)\frac{1}{v^2+(u-y)^2}+\big(-\frac{1}{2\lambda}\big)\frac{v(-2v)}{\big(v^2+(u-y)^2\big)^2}\Big)\,e^{-2\lambda v}\\
\frac{\frac{\partial}{\partial v}\Big(-\frac{1}{2\lambda}\frac{v}{v^2+(u-y)^2}e^{-2\lambda v}\Big)}{\frac{v}{v^2+(u-y)^2}e^{-2\lambda v}}&=1+\big(-\frac{1}{2\lambda}\frac{1}{v}\big)+\frac{1}{\lambda}\frac{v}{v^2+(u-y)^2}\\
&\in\Big(1-\frac{1}{2\lambda v},1+\frac{1}{\lambda v}\Big)\subseteq\Big(\frac{1}{2},2\Big),\cond{v\geqslant\frac{1}{\lambda}}.\qedhere\\
\endaligned
\]
\endproof

\begin{Cor}\label{lambda>0int} For any $r$ such that $0<r\leqslant e^{-\frac{1}{\lambda}}$,
\[
\aligned
\frac{1}{r^2}\int_{v>-\log r}H_\alpha(u+iv)|\!|\psi_\alpha'|\!|^2dv&\approx H_\alpha\big(u+(-\log r) i\big),\cond{0<|\alpha|<r^{1-\lambda}}\\
\frac{1}{r^2}\int_{v>\frac{\log|\alpha|-\log r}{\lambda}}H_\alpha(u+iv)|\!|\psi_\alpha'|\!|^2dv&\approx H_\alpha\big(u+(\textstyle\frac{\log|\alpha|-\log r}{\lambda}) i\big),\cond{0<|\alpha|<r^{1-\lambda}}\\
\frac{1}{r^2}\int_{v>\frac{\log|\alpha|-\log r}{\lambda}}H_\alpha(u+iv)|\!|\psi_\alpha'|\!|^2dv&\approx H_\alpha\big(u+(\textstyle\frac{-\log r}{\lambda}) i\big),\cond{0<|\alpha|<r^{1-\lambda}}
\endaligned
\]
\end{Cor}

\begin{figure}[H]
\begin{center}
   \includegraphics[width=0.4\linewidth]{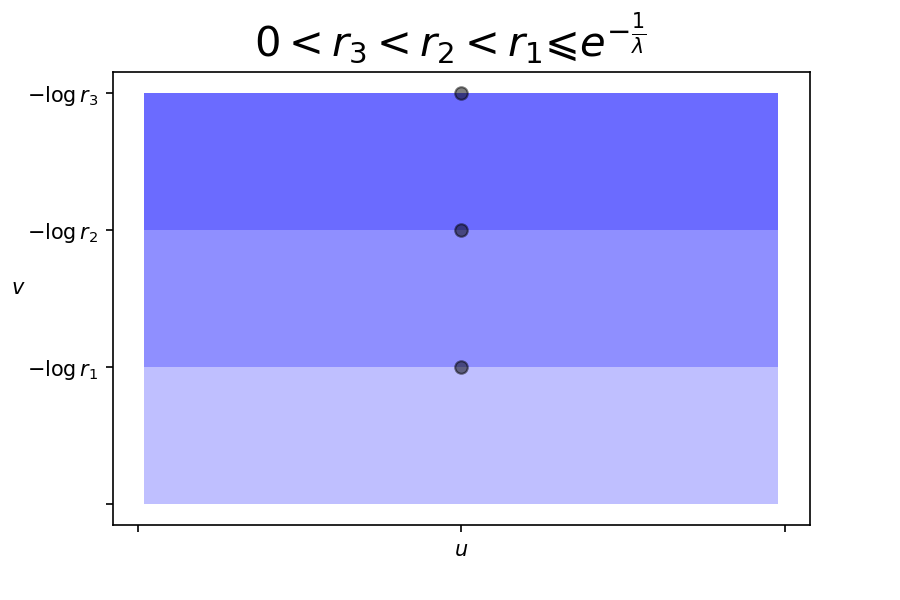}
   \caption{$\frac{1}{r^2}$ (The integration over $v>-\log r$) $\approx$ (The value at the boundary line $v=\log r$)}
   \end{center}
\end{figure}

\proof When $0<r\leqslant e^{-\frac{1}{\lambda}}$, which implies $-\log r\geqslant\frac{1}{\lambda}$, by Lemma \ref{approx-poission-lambda>0}
\[
\aligned
\int_{v>-\log r}H_\alpha(u+iv)|\!|\psi_\alpha'|\!|^2dv&=\frac{1}{\pi}\int_{v>-\log r}\int_{y\in\br}H_\alpha(y)\frac{v}{v^2+(u-y)^2}2\,(e^{-2v}+\lambda^2\,|\alpha|^2\,e^{-2\lambda v})dy\,dv\\
&\approx \frac{1}{\pi}\int_{y\in\br}H_\alpha(y)\Big\{\int_{v>-\log r}\frac{\partial}{\partial v}\Big(\frac{v}{v^2+(u-y)^2}\big(-e^{-2v}-\lambda\,|\alpha|^2\,e^{-2\lambda v}\big)\Big)dv\Big\}\,dy\\
&=\frac{1}{\pi}\int_{y\in\br}H_\alpha(y)\frac{-\log r}{(-\log r)^2+(u-y)^2}\big(r^2+\lambda\,|\alpha|^2\,r^{2\lambda}\big)dy\\
&=H_\alpha\big(u+(-\log r) i\big)\big(r^2+\lambda\,|\alpha|^2\,r^{2\lambda}\big)\\
&\approx r^2\,H_\alpha\big(u+(-\log r) i\big).
\endaligned
\]

For the same reason when $r^{1-\lambda}\leqslant |\alpha|<1$ which implies $\frac{\log|\alpha|-\log r}{\lambda}\geqslant-\log r\geqslant \frac{1}{\lambda}$
\[
\aligned
\int_{v>\frac{\log|\alpha|-\log r}{\lambda}}H_\alpha(u+iv)|\!|\psi_\alpha'|\!|^2dv
&\approx H_\alpha\big(u+(\textstyle\frac{\log|\alpha|-\log r}{\lambda}) i\big)\,\big(|\alpha|^{-\frac{2}{\lambda}}\,r^{\frac{2}{\lambda}}+\lambda\,r^2)\\
&\approx r^2\,H_\alpha\big(u+(\textstyle\frac{\log|\alpha|-\log r}{\lambda}) i\big).
\endaligned
\]

Last, when $|\alpha|\geqslant 1$ one has $\frac{-\log r}{\lambda}\geqslant-\log r\geqslant \frac{1}{\lambda}$ and
\[
\aligned
\int_{v>\frac{-\log r}{\lambda}}H_\alpha(u+iv)|\!|\psi_\alpha'|\!|^2dv
&\approx H_\alpha\big(u+(\textstyle\frac{-\log r}{\lambda}) i\big)\,\big(|\alpha|^{-\frac{2}{\lambda}}\,r^{\frac{2}{\lambda}}+\lambda\,r^2)\\
&\approx r^2\,H_\alpha\big(u+(\textstyle\frac{-\log r}{\lambda}) i\big).\qedhere
\endaligned
\]
\endproof

Thus
\[
\aligned
\cl(T,0)&\approx \lim\limits_{r\rw0+}\Big\{
\int_{|\alpha|<r^{1-\lambda}}\int_{u=0}^{2\pi}H_\alpha\big(u+(-\log r) i\big)\,du\,d\mu(\alpha)\\
& \ \ \ \ \ \ \ \ \ \ +\int_{r^{1-\lambda}\leqslant|\alpha|<1}\int_{u=0}^{2\pi}H_\alpha\big(u+(\textstyle\frac{\log|\alpha|-\log r}{\lambda}) i\big)\,du\,d\mu(\alpha)\\
& \ \ \ \ \ \ \ \ \ \ +\int_{|\alpha|\geqslant 1}\int_{u=0}^{2\pi}H_\alpha\big(u+(\textstyle\frac{-\log r}{\lambda}) i\big)\,du\,d\mu(\alpha)\Big\},
\endaligned
\]
by the inequalities (\ref{ineq1}) and (\ref{ineq2}) in the previous subsection. All terms are positive, so the order of taking limit and integration can change
\[
\aligned
\cl(T,0)&\approx\lim\limits_{v\rw+\infty}\int_{\alpha\in\bc^*}\int_{u=0}^{2\pi}H_\alpha\big(u+iv\big)\,du\,d\mu(\alpha)\\
&=\lim\limits_{k\rw+\infty}\int_{\alpha\in\bc^*}\int_{u=0}^{2\pi}\int_{y\in\br}H_\alpha(y)\frac{2k\pi}{(2k\pi)^2+(u-y)^2}\,dy\,du\,d\mu(\alpha),
\endaligned
\]

Fix some $k\in\bz$, $k\geqslant 2$. Define intervals $I_N$ for all $N\in\bz$ as follows
\begin{align*}
I_0&=[-2k\pi+2\pi,2k\pi),\\
I_N&=
\left\{
\begin{aligned}
&\big[2k\pi N,2k\pi(N+1)\big),  &\cond{N>0}\\
&\big[2k\pi (N-1)+2\pi,2k\pi N+2\pi\big),  &\cond{N<0}
\end{aligned}
\right.
\end{align*}
Thus $\br=\bigcup\limits_{N\in\bz}I_N$ is a disjoint union.
\begin{Lem} For any $u\in(0,2\pi)$ one has
\[
\frac{2k\pi}{(2k\pi)^2+(u-y)^2}\geqslant \frac{1}{1+(N+1)^2}\frac{1}{2k\pi},\cond{y\in I_N}
\]
\end{Lem}
\proof
Elementary.
\endproof

Thus
\[
\aligned
\cl(T,0)&\approx\lim\limits_{k\rw+\infty}\sum\limits_{N\in\bz}\int_{\alpha\in\bc^*}\int_{u=0}^{2\pi}\int_{y\in I_N}H_\alpha(y)\frac{2k\pi}{(2k\pi)^2+(u-y)^2}\,dy\,du\,d\mu(\alpha)\\
&\geqslant\lim\limits_{k\rw+\infty}\sum\limits_{N\in\bz}\int_{\alpha\in\bc^*}\int_{y\in I_N}\int_{u=0}^{2\pi}H_\alpha(y)\frac{1}{1+(N+1)^2}\frac{1}{2k\pi}\,du\,dy\,d\mu(\alpha)\\
&=\lim\limits_{k\rw+\infty}\sum\limits_{N\in\bz}\int_{\alpha\in\bc^*}\int_{y\in I_N}H_\alpha(y)\frac{1}{1+(N+1)^2}\frac{1}{k}\,dy\,d\mu(\alpha).
\endaligned
\]

By Lemma \ref{Lem:Mass}
\[
\aligned
\int_{\alpha\in\bc^*}\int_{y\in I_0}H_\alpha(y)\,dy\,d\mu(\alpha)&=(2k-1)\,|\!|T|\!|_{\bd^2}\\
&\geqslant k\,|\!|T|\!|_{\bd^2},\\
\int_{\alpha\in\bc^*}\int_{y\in I_N}H_\alpha(y)\,dy\,d\mu(\alpha)&=k\,|\!|T|\!|_{\bd^2},\cond{N\neq 0}.
\endaligned
\]
Thus
\[
\cl(T,0)\gtrsim\lim\limits_{k\rw+\infty}\sum\limits_{N\in\bz}\frac{1}{1+(N+1)^2}\,|\!|T|\!|_{\bd^2}\approx |\!|T|\!|_{\bd^2}
\]
is nonzero.

\section{Case $\lambda<0$. Periodic current, including all currents when $\lambda\in\bq_{<0}$}
For any $\alpha\in\bc^*$ fixed, the leaf $L_\alpha$ is contained in a real 3-dimensional analytic Levi-flat CR manifold $|w|=|\alpha|\,|z|^{\lambda}$, which can be viewed as a curve in $|z|,|w|$ coordinates. The norms $|z|$ and $|w|$ depends only on $v$. No leaf $L_\alpha$ tends to the singularity $(0,0)$. For $r$ sufficiently small, the leaf $L_\alpha$ is outside of $r\,\bd^2$.

\begin{figure}[H]
\begin{center}
   \includegraphics[width=0.45\linewidth]{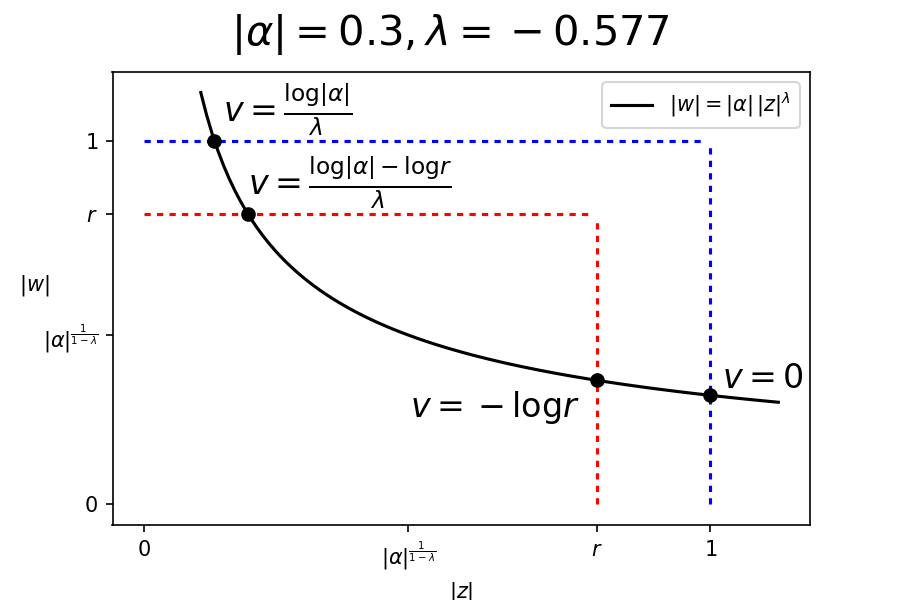}
   \caption{Case $\lambda<0$}
\end{center}
\end{figure}

Like the case $\lambda>0$, when one fixes $|z|=r$ for some $r\in(0,1)$, then $|w|=|\alpha|\,|z|^\lambda$ is uniquely determined and the real 2-dimensional leaf $L_\alpha$ becomes a real 1-dimensional curve $L_{\alpha,r}:=L_\alpha\cap \bt^2_r$ on the torus $\bt^2_r:=\{(z,w)\in\bd^2~|~|z|=r,|w|=|\alpha|\,r^\lambda\}$. It is a closed curve if $\lambda\in\bq$, a dense curve on $\bt^2_r$ if $\lambda\notin\bq$.

Let $T$ be a harmonic current directed by $\cf$. Then $T|_{P_\alpha}$ has the form $h_\alpha(z,w)[P_\alpha]$. Let $H_\alpha:=h_\alpha\circ\psi_\alpha(u+iv)$. It is a positive harmonic function for $\mu$-almost all $\alpha\in\bd^*$ defined on a neighborhood of a horizontal strip $\{(u,v)\in\br^2~|~0<v<\frac{\log|\alpha|}{\lambda}\}$. 

Like what happened in the case $\lambda>0$, one only calculates the mass on an open subset $U:=\{(z,w)\in\bd^2~|~z\notin\br_{\geqslant0},w\neq 0\}$. For each $\alpha\in\bd^*$ one normalizes $H_\alpha$ by setting $H_\alpha(0)=1$ to fix the expression of $T:=\int h_\alpha [P_\alpha] d\mu(\alpha)$. Similar to Lemma \ref{Lem:Mass}, for each $k_0\in\bz$ fixed,

\[
\aligned
|\!|T|\!|_{\bd^2}&=\int_{0<|\alpha|<1}\int_{v=0}^{\frac{\log|\alpha|}{\lambda}}\int_{u=2\,k_0\,\pi}^{2\,k_0\,\pi+2\,\pi}H_\alpha(u+iv)\,2\,(e^{-2v}+\lambda^2\,|\alpha|^2\,e^{-2\lambda\,v})\,du\,dv\,d\mu(\alpha)\\
\cl(T,0)&=\lim\limits_{r\rw 0+}\frac{1}{r^2}|\!|T|\!|_{r\bd^2}\\
&=\lim\limits_{r\rw 0+}\frac{1}{r^2}\int_{0<|\alpha|<r^{1-\lambda}}\int_{v=-\log r}^{\frac{\log|\alpha|-\log r}{\lambda}}\int_{u=2\,k_0\,\pi}^{2\,k_0\,\pi+2\,\pi}H_\alpha(u+iv)\,2\,(e^{-2v}+\lambda^2\,|\alpha|^2\,e^{-2\lambda\,v})\,du\,dv\,d\mu(\alpha)
& \ \ \ \ \ \ \ \ \ \ \ \ \ \ \ \ \ \ \ \ \ \ \ \ \ \ \ \ \ \ \ \ \ \ \ \ \ \ \ \ \ \ \ \ \ \ \ \ \ \ \ \ \ \ \ \ \ \ \ \ \ \ \ \ \ \ \ \ \ \ \ \ \ \ \ \ \ \ \ \ \ \ \ \ \ \ \ \ \ \ \ \ \ \ \ \ \ \ \ \ \ \      dy\,dx\,d\mu(\alpha)
\endaligned
\]

Now assume $T$ is periodic, we treat Theorem~\ref{thm:negative-periodic}. Suppose that there exists some $b\in\bz_{\leqslant 1}$ such that $H_{\alpha}(u+iv)=H_{\alpha}(u+2\pi b+iv)$ for all $\alpha\in\bd^*$ and all $(u,v)$ in a neighborhood of the strip $\{(u+iv)\in\bc~|~u\in\br,v\in[0,\frac{\log|\alpha|}{\lambda}]\}$. Like Lemma \ref{lem:periodic}, one proves

\begin{Lem} Let $F(u,v)$ be a positive harmonic function on a neighborhood of the horizontal strip $\{(u+iv)\in\bc~|~u\in\br,v\in[0,C]\}$ for some $C>0$. Suppose $F(u,v)=F(u+2\pi b,v)$ on this strip. Then
\[
F(u,v)=\sum\limits_{k\in\bz,k\neq 0}\big(a_k\,e^{\frac{kv}{b}}\cos(\tfrac{ku}{b})+b_k\,e^{\frac{kv}{b}}\sin(\tfrac{ku}{b})\big)+a_0\,(1-C^{-1}\,v)+b_0\,v,
\]
for some $a_k,b_k\in\br$ with $a_0\geqslant0$ and $b_0\geqslant0$.
\end{Lem}
\proof The proof is almost the same as that of Lemma~\ref{lem:periodic}. Using Fourier series and calculating Laplacian, one concludes that
\[
F(u,v)=\sum\limits_{k\in\bz,k\neq 0}\big(a_k\,e^{\frac{kv}{b}}\cos(\tfrac{ku}{b})+b_k\,e^{\frac{kv}{b}}\sin(\tfrac{ku}{b})\big)+p+q\,v,
\]
for some $a_k,b_k,p,q\in\br$. For any $v\in[0,C]$, $F(u,v)\geqslant 0$ implies
\[
\int_{u=0}^{2\pi b}F(u,v)du=2\pi b\, (p+q\,v)\geqslant0.
\]
Thus $p\geqslant0$ and $q\geqslant-C^{-1}\,p$. One may write $p+q\,v=p\,(1-C^{-1}\,v)+(q+C^{-1}\,p)\,v$ with $p=:a_0\geqslant 0$ and $q+C^{-1}\,p=:b_0\geqslant 0$.\endproof

For periodic currents one may assume
\begin{align}
\label{periodic-H-negative}
H_\alpha(u+iv)=\sum\limits_{k\in\bz,k\neq 0}\big(a_k(\alpha)\,e^{\frac{kv}{b}}\cos(\tfrac{ku}{b})+b_k(\alpha)\,e^{\tfrac{kv}{b}}\sin(\tfrac{ku}{b})\big)+a_0(\alpha)\,(1-\textstyle\frac{\lambda}{\log|\alpha|}\,v)+b_0(\alpha)\,v,
\end{align}
for some $a_k(\alpha),b_k(\alpha)\in\br$ with $a_0(\alpha)\geqslant0$ and $b_0(\alpha)\geqslant0$. According to Lemma~\ref{Lem:Mass}, for any $k_0\in\bz$, use the jacobian \eqref{jacobian}
\[
|\!|T|\!|_{\bd^2}=\int_{0<|\alpha|<1}\int_{v=0}^{\frac{\log|\alpha|}{\lambda}}\int_{u=2k_0\pi}^{2k_0\pi+2\pi}H_\alpha(u+iv)\,2\,(e^{-2v}+\lambda^2\,|\alpha|^2\,e^{-2\lambda\,v})\,du\,dv\,d\mu(\alpha).
\]
Next, using $0=\int_0^{2\pi b}\cos(\frac{ku}{b})du$ for $k\neq 0$ and the same for $\sin(\frac{ku}{b})$, let us calculate the average among $k_0=0,1,\dots,b-1$ for the mass
\[
\aligned
|\!|T|\!|_{\bd^2}&=\frac{1}{b}\int_{0<|\alpha|<1}\int_{v=0}^{\frac{\log|\alpha|}{\lambda}}\int_{u=0}^{2\pi b}H_\alpha(u+iv)\,2\,(e^{-2v}+\lambda^2\,|\alpha|^2\,e^{-2\lambda\,v})\,du\,dv\,d\mu(\alpha)\\
&=\frac{2\pi b}{b}\int_{0<|\alpha|<1}\int_{v=0}^{\frac{\log|\alpha|}{\lambda}}\left(a_0(\alpha)\,(1-\textstyle\frac{\lambda}{\log|\alpha|}\,v)+b_0(\alpha)\,v\right)\,2\,(e^{-2v}+\lambda^2\,|\alpha|^2\,e^{-2\lambda\,v})\,dv\,d\mu(\alpha),
\endaligned
\]
and for the Lelong number
\[
\aligned
\cl(T,0)&=\lim\limits_{r\rw 0+}\frac{1}{r^2}|\!|T|\!|_{r\bd^2}\\
&=\lim\limits_{r\rw 0+}\frac{1}{b\,r^2}\int_{0<|\alpha|<r^{1-\lambda}}\int_{v=-\log r}^{\frac{\log|\alpha|-\log r}{\lambda}}\int_{u=0}^{2\pi b}H_\alpha(u+iv)\,2\,(e^{-2v}+\lambda^2\,|\alpha|^2\,e^{-2\lambda\,v})\,du\,dv\,d\mu(\alpha)\\
&=\lim\limits_{r\rw 0+}\frac{2\pi b}{b\,r^2}\int_{0<|\alpha|<r^{1-\lambda}}\int_{v=-\log r}^{\frac{\log|\alpha|-\log r}{\lambda}}\left(a_0(\alpha)\,(1-\textstyle\frac{\lambda}{\log|\alpha|}\,v)+b_0(\alpha)\,v\right)\,2\,(e^{-2v}+\lambda^2\,|\alpha|^2\,e^{-2\lambda\,v})\,dv\,d\mu(\alpha).
\endaligned
\]

Introduce the two functions of $r\in(0,1]$ given by elementary integrals
\[
\aligned
I_a(r)&:=\frac{1}{r^2}\int_{v=-\log r}^{\frac{\log|\alpha|-\log r}{\lambda}}2\,(1-\textstyle\frac{\lambda}{\log|\alpha|}\,v)\,(e^{-2v}+\lambda^2\,|\alpha|^2\,e^{-2\lambda\,v})\,dv,\\
&=1+\lambda\,|\alpha|^2 \,  r^{2 \lambda -2}+\frac{1}{2\log|\alpha|}\left(-2|\alpha|^{-\frac{2}{\lambda} } r^{\frac{2}{\lambda }-2} \log (r)+\lambda  |\alpha|^{-\frac{2}{\lambda} } r^{\frac{2}{\lambda }-2}+2\,\lambda ^2\,|\alpha|^2 \, r^{2 \lambda -2} \log (r)-\lambda\,|\alpha|^2 \,  r^{2 \lambda -2}\right),\\
I_b(r)&:=\frac{1}{r^2}\int_{v=-\log r}^{\frac{\log|\alpha|-\log r}{\lambda}}2\,v\,(e^{-2v}+\lambda^2\,|\alpha|^2\,e^{-2\lambda\,v})\,dv\\
&=\frac{1}{2} \left(-\frac{|\alpha|^{-\frac{2}{\lambda} } r^{\frac{2}{\lambda }-2} (\lambda+2 \log |\alpha| -2 \log (r))}{\lambda }+|\alpha|^2 r^{2 \lambda -2} (1-2 \lambda  \log (r))-2 \log |\alpha|\right),
\endaligned
\]
to describe the contributions from $a_0(\alpha)$ part and from $b_0(\alpha)$ part.

   \begin{figure}[H]
   \begin{center}
   \includegraphics[width=0.3\linewidth]{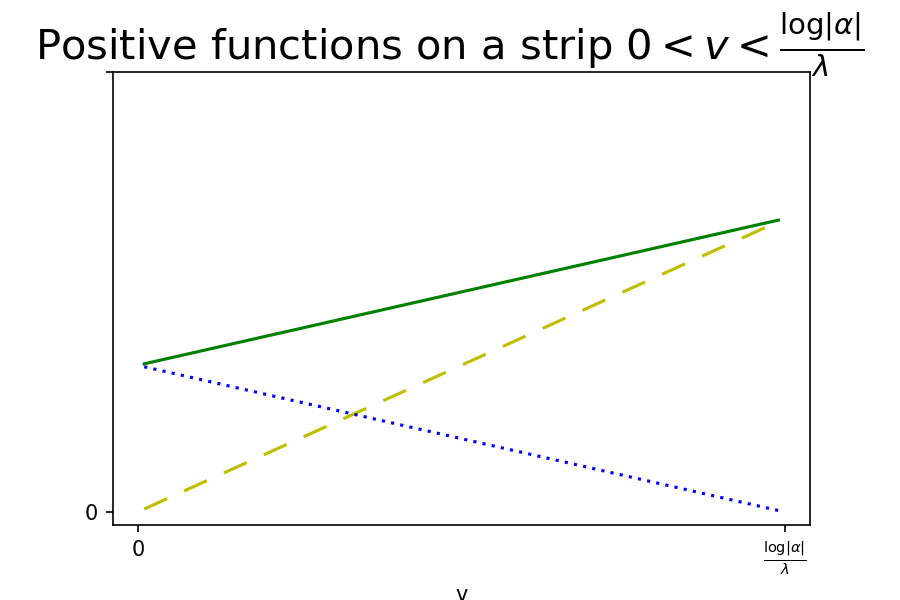}
   \caption{The green = the blue dotted (gives $I_a(r)$) + the yellow dashed (gives $I_b(r)$)}
   \end{center}
   \end{figure}

Then we can express
\[
\aligned
|\!|T|\!|_{\bd^2}&=2\,\pi\,\int_{0<|\alpha|<1}\Big(a_0(\alpha)\,I_a(1)+b_0(\alpha)\,I_b(1)\Big)\,d\mu(\alpha),\\
\cl(T,0)&=2\,\pi\,\lim\limits_{r\rw 0+}\int_{0<|\alpha|<r^{1-\lambda}}\Big(a_0(\alpha)\,I_a(r)+b_0(\alpha)\,I_b(r)\Big)\,d\mu(\alpha).
\endaligned
\]
Observe that
\[
\aligned
I_a(1)&=1+\lambda\,|\alpha|^2 +\frac{\lambda  \left(|\alpha|^{-\frac{2}{\lambda} }-|\alpha|^2\right)}{2 \log |\alpha|},\\
I_b(1)&=\frac{1}{2} \left(-\frac{|\alpha|^{-\frac{2}{\lambda} } (\lambda+2 \log|\alpha| )}{\lambda }+|\alpha|^2-2 \log|\alpha|\right).
\endaligned
\]

Fix any $\alpha\in\bd^*$, by definition $r^2I_a(r)$ and $r^2I_b(r)$ are increasing for $r\in(0,1]$, since the interval of integration $(-\log r,\frac{\log|\alpha|-\log r}{\lambda})$ is expanding and the function integrated is positive. In particular, for any $r\in(0,1]$,
\[
I_a(r)\leqslant r^{-2}\,I_a(1), \ \ \ \ I_b(r)\leqslant r^{-2}\,I_b(1).
\]

It is more subtle to talk about monotonicity of $I_a(r)$ and $I_b(r)$. We expect upper bounds of $I_a(r)/I_a(1)$ and $I_b(r)/I_b(1)$ for $r\in(0,1]$ which are independent of $\alpha$, i.e. depend only on $\lambda$.

\begin{Lem} For any $r\in(0,1)$ and any $\alpha\in\bc$ with $0<|\alpha|<r^{1-\lambda}<1$, one has
\[
0<I_a(r)<I_a(1).
\]
\end{Lem}
\proof
A differentiation gives
\[
\aligned
\frac{d}{dr}I_a(r)&=\underbrace{\frac{|\alpha|^{-\frac{2}{\lambda} }}{\lambda\,r^3\,\log |\alpha|}}_{>0}
\Big(\lambda ^2 \big(|\alpha|^{2+\frac{2}{\lambda }} r^{2 \lambda }-r^{\frac{2}{\lambda} }\big)-2 (1-\lambda) \big(\lambda ^3 |\alpha|^{2+\frac{2}{\lambda }} r^{2 \lambda }+r^{\frac{2}{\lambda} }\big)\log (r) 
\\
& \ \ \ \ \ \ \ \ \ \ \ \ \ \ \ \ \ \  -2 (1-\lambda) \lambda ^2 |\alpha|^{2+\frac{2}{\lambda }} r^{2 \lambda }\log|\alpha| \Big).
\endaligned
\]
It suffices to show that $\frac{d}{dr}I_a(r)>0$ when $r\in(0,1)$ and $0<|\alpha|<r^{1-\lambda}$.

Introduce the new variable $t:=\frac{|\alpha|}{r^{1-\lambda}}\in(0,1)$. In the big parentheses, replace $|\alpha|$ by $t\,r^{1-\lambda}$ and $\log|\alpha|$ by $\log(t)+(1-\lambda)\log(r)$
\[
\aligned
\frac{d}{dr}I_a(r)&=\underbrace{\frac{|\alpha|^{-\frac{2}{\lambda} }\,r^{\frac{2}{\lambda}}}{\lambda  r^3 \log |\alpha|}}_{>0}
\Big(
\lambda^2(t^{2+\frac{2}{\lambda}}-1)
-2\,(1-\lambda)\,(t^{2+\frac{2}{\lambda}}+1)\log (r)\\
& \ \ \ \ \ \ \ \ \ \ \ \ \ \ \ \ \ \ \underbrace{-2\,(1-\lambda)\,\lambda^2\,t^{2+\frac{2}{\lambda}}\log (t)}_{>0}
\Big)\\
&>\frac{|\alpha|^{-\frac{2}{\lambda} }\,r^{\frac{2}{\lambda}}}{\lambda  r^3 \log |\alpha|}
\Big(
\lambda^2\,\underbrace{(t^{2+\frac{2}{\lambda}}-1)}_{\geqslant0}
\underbrace{-2\,(1-\lambda)\,(t^{2+\frac{2}{\lambda}}+1)\log (r)}_{>0}\Big)>0,
\endaligned
\]
since $\lambda\in[-1,0)$ implies $t^{2+\frac{2}{\lambda}}\geqslant 1$.\endproof

It is not true that $I_b(r)$ is increasing on $(0,1]$, but on a smaller half-neighborhood of $0$, independent of $\alpha$, it is increasing. This suffice to give an upper bound of $I_b(r)/I_b(1)$.

\begin{Lem} For any $r\in(0,e^{\frac{1}{2\,\lambda\,(1-\lambda)}})$ and any $\alpha\in\bc$ with $0<|\alpha|<r^{1-\lambda}<1$, one has
\[
0<I_b(r)<I_b(e^{\frac{1}{2\lambda(1-\lambda)}})\leqslant e^{\frac{1}{-\lambda\,(1-\lambda)}}\,I_b(1).
\]
\end{Lem}
\proof
A differentiation gives
\[
\aligned
\frac{d}{dr}I_b(r)&=\underbrace{\frac{|\alpha|^{-\frac{2}{\lambda} }}{\lambda ^2 r^3}}_{>0}
\Big(-\lambda ^2 \big(|\alpha|^{2+\frac{2}{\lambda }} r^{2 \lambda }-r^{\frac{2}{\lambda}}\big) 
+2 (1-\lambda)  \big(\lambda ^3 |\alpha|^{2+\frac{2}{\lambda }} r^{2 \lambda }+r^{\frac{2}{\lambda} }\big) \log (r)\\
& \ \ \ \ \ \ \ \ \ \ \ \ \ \ \ \ \ \  -2 (1-\lambda)  r^{\frac{2}{\lambda} } \log |\alpha|\Big)
\endaligned
\]
It suffices to show that $\frac{d}{dr}I_b(r)>0$ when $0<r<e^{\frac{1}{2\,\lambda\,(1-\lambda)}}$ and $0<|\alpha|<r^{1-\lambda}$.

Again, introduce the variable $t:=\frac{|\alpha|}{r^{1-\lambda}}\in(0,1)$ and replace $\alpha$ and $\log|\alpha|$ in the parentheses
\begin{align*}
\frac{d}{dr}I_b(r)&=\underbrace{\frac{|\alpha|^{-\frac{2}{\lambda} }\,r^{\frac{2}{\lambda}}}{\lambda^2\,r^3}}_{>0}
\Big(
-\lambda^2\,(t^{2+\frac{2}{\lambda}}-1)+2\,\lambda\,(1-\lambda)\,(\lambda^2\,t^{2+\frac{2}{\lambda}}+1)\,\log(r)
\\
& \ \ \ \ \ \ \ \ \ \ \ \ \ \ \ \ \ \  \underbrace{-2 (1-\lambda ) \log(t) }_{>0}\Big)\\
&>\frac{|\alpha|^{-\frac{2}{\lambda} }\,r^{\frac{2}{\lambda}}}{\lambda^2\,r^3}
\Big(
-\lambda^2\,(t^{2+\frac{2}{\lambda}}-1)+\underbrace{2\,\lambda\,(1-\lambda)\,(\lambda^2\,t^{2+\frac{2}{\lambda}}+1)}_{<0}\,\underbrace{\log(r)}_{<\frac{1}{2\lambda(1-\lambda)}<0}\Big)\\
&>\frac{|\alpha|^{-\frac{2}{\lambda} }\,r^{\frac{2}{\lambda}}}{\lambda^2\,r^3}
\Big(
-\lambda^2\,(t^{2+\frac{2}{\lambda}}-1)+\lambda^2\,t^{2+\frac{2}{\lambda}}+1\Big)=\frac{|\alpha|^{-\frac{2}{\lambda} }\,r^{\frac{2}{\lambda}}}{\lambda^2\,r^3}
\big(
\lambda^2\,+1\big)>0.\qedhere
\end{align*}
\endproof

\proof[End of proof of Theorem~\ref{thm:negative-periodic}]
From what precedes, the Lelong number is zero
\[
\aligned
\cl(T,0)&=2\,\pi\,\lim\limits_{r<e^{\frac{1}{2\lambda(1-\lambda)}},r\rw 0+}\int_{0<|\alpha|<r^{1-\lambda}}\big(a_0(\alpha)\,I_a(r)+b_0(\alpha)\,I_b(r)\big)\,d\mu(\alpha)\\
&<2\,\pi\,\lim\limits_{r\rw 0+}\int_{0<|\alpha|<r^{1-\lambda}}\big(a_0(\alpha)\,I_a(1)+b_0(\alpha)\,e^{\frac{1}{-2\,\lambda\,(1-\lambda)}}\,I_b(1)\big)\,d\mu(\alpha)\\
&\approx2\,\pi\,\lim\limits_{r\rw 0+}\int_{0<|\alpha|<r^{1-\lambda}}\big(a_0(\alpha)\,I_a(1)+b_0(\alpha)\,I_b(1)\big)\,d\mu(\alpha)=0,
\endaligned
\]
since $|\!|T|\!|_{\bd^2}=2\,\pi\,\int_{0<|\alpha|<1}\big(a_0(\alpha)\,I_a(1)+b_0(\alpha)\,I_b(1)\big)\,d\mu(\alpha)$ is finite.\endproof

%%%%%%%%%%%%%%%%%%%%%%%%%%%%
\setlength\parindent{0em}
{\scriptsize Zhangchi Chen, Universit\'e Paris-Saclay, CNRS, Laboratoire de math\'ematiques d'Orsay, 91405, Orsay, France}\\
{\bf\scriptsize zhangchi.chen@universite-paris-saclay.fr}, {\bf\scriptsize https://www.imo.universite-paris-saclay.fr/$\sim$chen/}
%%%%%%%%%%%%%%%%%%%%%%%%%%%%
\end{document}